\documentclass{book}

\usepackage{amsmath, amsfonts, amssymb, graphicx, amsthm}

\newtheorem{theorem}{Theorem}[section]
\newtheorem{lemma}[theorem]{Lemma}
\newtheorem{proposition}[theorem]{Proposition}

\newtheorem{definition}[theorem]{Definition}

\newtheorem{remark}[theorem]{Remark}
\newtheorem{example}[theorem]{Example}

\begin{document}

\title{A short course on positive solutions of systems of ODEs via fixed point index}
\author{Gennaro Infante\\
Dipartimento di Matematica e Informatica\\ Universit\`{a} della
Calabria\\ 87036 Arcavacata di Rende, Cosenza, Italy\\ www.mat.unical.it/$\sim$infante\\ gennaro.infante@unical.it}%
\date{}%
\maketitle

\tableofcontents

\chapter*{On this short course}\addcontentsline{toc}{chapter}{On this short course}
These notes were used for a mini-course delivered during the Workshop ``Differential Equations and Applications'' held at the Campus of Ourense of the University of Vigo (Spain) in June 2016. A shorter version of this course was previously given to doctoral students in Santiago De Compostela (Spain) in May 2013 and September 2014 and in Ruse (Bulgaria) in September 2013, within the framework of the Erasmus+ program.

This short course is meant for doctoral students and young researchers interested in the application of topological methods to the solvability of differential and integral equations. 

The notes are organized as follows. In Chapter~\ref{claskras} we discuss the existence of one positive solution of a model problem (a simple second order ODE with Dirichlet boundary conditions)
via the classical Krasnosel'ski\u\i{} fixed-point theorem. In Chapter~\ref{fpin} we illustrate how to prove existence and multiplicity results for the model problem via the fixed point index theory for compact maps. We also provide, by means of elementary arguments, some non-existence results. In Chapter~\ref{secsyst} we show how the approach developed for one equation can be tailored in order to deal with the existence and multiplicity of non-negative solutions for systems of ODEs subject to local boundary conditions. In Chapter~\ref{mgenbcs} we discuss the case of more general boundary conditions, focusing on a three-point problem and on nonlinear boundary conditions. In Chapter~\ref{radsolpdes} we briefly illustrate how to adapt the theory in order to deal with the existence of radial solutions of some systems of elliptic PDEs subject to local and nonlocal boundary conditions in the case of annular or exterior domains.

\chapter[The Krasnosel'ski\u\i{} fixed point theorem]{The classical Krasnosel'ski\u\i{} fixed point theorem}\label{claskras}
A classical problem is to investigate the existence of positive solutions for the second order differential equation 
\begin{equation}\label{EQ-dir}
u''(t)+f(u(t))=0,\ t \in (0,1),\\
\end{equation}
subject to \emph{Dirichlet} boundary conditions (BCs)
\begin{equation}\label{BC-dir}
u(0)=u(1)=0,
\end{equation}
where $f$ is a continuous function.
 
One motivation is that this problem often occurs when studying the existence of radial solutions in $\mathbb{R}^n,\; n\ge 2$, for the boundary value problem (BVP)
\begin{equation*}\label{eq6.1}
\triangle v+ f(v) = 0,\; x\in \mathbb{R}^n,\; \ |x|\in
[R_{1},R_{2}],
\end{equation*}
\par\noindent
with
\begin{equation*}
v=0 \quad \text{for}  \quad
|x|=R_{1}\quad \text{and}\quad  |x|=R_{2},
\end{equation*}
where $0<R_{1}<R_{2}<\infty$.

Several methods have been used to study the BVP~\eqref{EQ-dir}-\eqref{BC-dir}, for example upper and lower solutions, variational methods and shooting methods. 

We begin by considering a well-known tool, the fixed point theorem of \emph{Krasnosel'ski\u\i}, sometimes called ``the cone compression-expansion Theorem''.
\begin{definition}
A cone $K$ in a Banach space $X$ is a closed convex set such that $\lambda x\in K$ for every $x \in K$ and for all $\lambda\geq 0$ and satisfying $K\cap (-K)=\{0\}$.
\end{definition}
\begin{example} Two examples of cones:
\begin{enumerate}
\item In $\mathbb{R}^2$, the set $\mathbb{R}^2_{+}:=\{(x,y)\in \mathbb{R}^2\,\text{such that}\, x\geq 0,\, y\geq 0\}$ is a cone.
\item In $C[0,1]$, the set 
$
P: = \bigl\{u\in C[0,1]:  u(t)\geq 0\bigr\} 
$
is a cone.
\end{enumerate}
\end{example}
\begin{theorem}[Krasnosel'ski\u\i{}-Guo, (1962; 1985)]\label{krasno}
Let $T:K \to K$ be a compact map\footnote{By compact we mean that $T$ is
continuous and $\overline {T(Q)}$ is compact for each bounded
subset $Q\subset K$.}. Assume that there exist two positive constants $r,R$ with $r\neq R$ such that 
\begin{align*}
 & \|Tu\| \leq \|u\| \; \text{for every} \;
u \in K \;\text{with}\; \|u\|=r, \\
 & \|Tu\| \geq \|u\| \; \text{for every} \;
u \in K \;\text{with}\; \|u\|=R.
\end{align*}
Then there exists  $u_0\in K$ such that  $Tu_0=u_0$ and 
$
\min\{r,R\}\leq \|u\| \leq \max\{r,R\}.
$
\end{theorem}
\begin{figure}[h]%\label{fig1}
\centering
\includegraphics[height=5.5cm,angle=0]{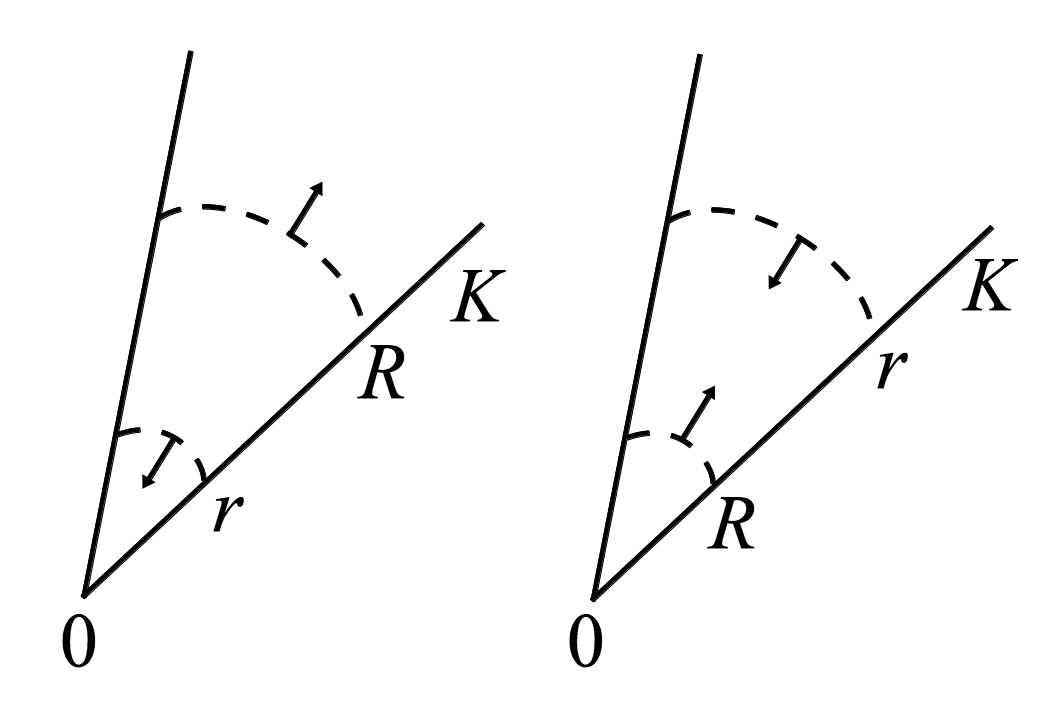}
\caption{An illustration of the cone compression-expansion Theorem.}
\end{figure}{}
We postpone the proof of the theorem, which follows from classical fixed point index theory for compact maps, and we focus on how to use it to study our problem:
the idea  is to rewrite our BVP as an integral equation in a suitable space. This is not too dissimilar to what happens when initial value problems are rewritten in the form of Volterra equations. In particular we would like to rewrite the BVP as a \emph{Hammerstein} integral equation
$$
u(t)=\int_{0}^{1}k(t,s)f(u(s))ds,
$$
where the function $k$ is said to be the \emph{kernel} of the integral equation (or the \emph{Green's function}\footnote{Named after the British mathematical physicist George Green (1793-1841).} of the problem). There are many ways of constructing the Green's function (for example by variation of parameters or by Laplace transforms), in our case we proceed as follows.

Consider the linear problem
\begin{align}\label{linBVP}
u''(t)+y(t)=0,\ u(0)=u(1)=0.
\end{align}
If we integrate $u''+y=0$ we obtain
$$u'(t)=u'(0)+\int_{0}^{t}-y(s)ds$$ and, integrating again, we get
$$u(t)=u(0)+tu'(0)-\int_{0}^{t}\int_{0}^{w}y(s)dsdw.$$
Using the Cauchy formula for iterated kernels\footnote{In general the formula reads as follows: $$\int_{0}^t\int_0^{s_1}\cdots\int_0^{s_{n-1}}f(s_{n})ds_{n}\cdots ds_2 ds_1=\frac{1}{(n-1)!}\int_0^t\left(t-s\right)^{n-1}f(s)ds.$$}, we obtain 
$$u(t)=u(0)+tu'(0)-\int_{0}^{t}(t-s)y(s)ds.$$
We now make use of the boundary conditions
\begin{align*}
u(t)&=tu'(0)-\int_{0}^{t}(t-s)y(s)ds,\\
u(1)&=u'(0)-\int_{0}^{1}(1-s)y(s)ds.
\end{align*}
This gives
$$
u(t)=\int_{0}^{1}t(1-s)y(s)ds-\int_{0}^{t}(t-s)y(s)ds.
$$
Now 
$$
\int_{0}^{t}t(1-s)y(s)ds-\int_{0}^{t}(t-s)y(s)ds=\int_{0}^{t}(t-st-t+s)y(s)ds=\int_{0}^{t}s(1-t)y(s)ds,
$$
which yields
$$
u(t)=\int_{0}^{1}k(t,s)y(s)ds,
$$
where
\begin{equation}\label{kDir}
k(t,s)=\begin{cases}
  s(1-t)& \text{if}\;s \leq t,\\
  t(1-s) &\text{if}\; s>t.
\end{cases}
\end{equation}

Once we have found the Green's function for \eqref{linBVP}, the integral equation associated to the BVP~\eqref{EQ-dir}-\eqref{BC-dir} is
given by
\begin{equation}\label{eqham}
u(t)=\int_{0}^{1}k(t,s)f(u(s))ds.
\end{equation}
If $f$ is continuous it can be proved that $u$ is a solution of the integral equation~\eqref{eqham} if and only if $u$ is a solution of the BVP~\eqref{EQ-dir}-\eqref{BC-dir}.

We want to study the solutions of the equation \eqref{eqham} as fixed points of the Hammerstein integral operator 
\begin{equation*}%\label{eqhamop}
Tu(t):=\int_{0}^{1}k(t,s)f(u(s))ds,
\end{equation*}
in a suitable space. Here we consider $C[0,1]$, endowed with the usual supremum norm $\|u\|:=\max_{t\in [0,1]}|u(t)|$, and we assume that 
\begin{itemize}
\item  $f:[0, +\infty) \to [0, +\infty)$ is continuous.
\end{itemize}
A natural setting would be to look for fixed points of the operator $T$ in the cone 
$$P= \bigl\{u\in C[0,1]:  u(t)\geq 0\bigr\}.$$
We will show the existence of positive solutions in a type of cone, introduced by D. Guo \cite{guo}, which is smaller than $P$, namely
$$
K = \{u\in P:  \displaystyle{\min_{t\in[a,b]}}u(t) \geq c \|u\|\},
$$
where $[a,b]\subseteq[0,1]$ and $c\in(0,1]$. 

In order to find the interval $[a,b]$ and the constant $c$, we look for upper and lower bounds for the Green's function $k$; in other words, we look for  
a continuous function $\Phi:[0,1]\to [0, +\infty)$ and a number $c\in(0,1]$ such that
\begin{align*}
  k(t,s)\leq \Phi(s) \text{ for }& t,s \in [0,1],\\
  c \Phi(s)\leq k(t,s) \text{ for }& t \in [a,b] \text{ and }s \in
  [0,1].
\end{align*}{}
Now
$$k(t,s)=\begin{cases}
  s(1-t)& \text{if}\;s \leq t,\\
  t(1-s) &\text{if}\; s>t,
\end{cases}
$$
therefore we have 
$$k(t,s)\leq s(1-s) \text{ for } t,s \in [0,1].$$
Now let $t\in [a,b]$ and $s \in [0,1]$. If $s\leq t$ we have $$k(t,s)\geq s(1-b)\geq (1-b)s(1-s),$$ and if $s>t$ we have  $$k(t,s)\geq a(1-s)\geq as(1-s).$$
Thus we may choose $\Phi(s)=s(1-s)$, $[a,b]=[1/4, 3/4]$ and $c=1/4$ and we can work in the cone
$$
K = \bigl\{u\in C[0,1]:  u\geq 0, \displaystyle{\min_{t\in[1/4,3/4]}}u(t) \geq \frac{1}{4} \|u\|\bigr\}.
$$
In order to apply Theorem~\ref{krasno} we need to show the following.
\begin{lemma}
The operator $T$ maps $K$ into $K$ and is compact.
\end{lemma}
\begin{proof}
We show that $T:K\to K$. Indeed, we have, for $t\in [0,1]$,
$$
|Tu(t)| = \int_0^1 k(t,s)f(u(s))\,ds\leq  \int_{0}^{1}\Phi(s)f(u(s))\,ds
$$
so that
$$ \|Tu\|\leq \int_{0}^{1}\Phi(s)f(u(s))\,ds.
$$
Also we have $$ \min_{t\in [1/4,3/4]}\{Tu(t)\} \geq  \frac{1}{4}
\int_{0}^{1}\Phi(s)f(u(s))\,ds. $$ Hence $Tu\in K$ for every
$u\in K$. 

The compactness of $T$ follows from the classical Ascoli-Arzel\`{a} Theorem.
\end{proof}

\begin{lemma}\label{lemcom}
Assume that there exists $\rho> 0$ such that
$
f^{0,\rho} \leq m,
$
where
\begin{equation*}
  f^{0,{\rho}}:=\sup\Bigl\{\frac{f(u)}{\rho} : 0 \leq u \leq {\rho}\Bigr\}\ \text{and}\
\frac{1}{m}:=\sup_{t\in [0,1]} \int_{0}^{1} k(t,s) \,ds.
\end{equation*}
Then $\|Tu\| \leq \|u\| $ for every $u \in K \;\text{with}\; \|u\|=\rho$.
\end{lemma}

\begin{proof}
Take $u\in K$ with $\|u\|=\rho$. Then for $t\in [0,1]$ we have
$$Tu(t)=\int_{0}^{1}k(t,s)f(u(s))ds\leq \rho  f^{0,{\rho}}\int_{0}^{1} k(t,s) \,ds\leq \rho=\|u\|.$$
\end{proof}
\begin{lemma}\label{lemexp}
Assume that there exists $\rho> 0$ such that
$f_{{\rho/4},{\rho }}
\geq M,
$
where
\begin{equation*}
  f_{{\rho}/{4},{\rho}}:=\inf\Bigl\{\frac{f(u)}{\rho} : \rho/4 \leq u \leq {\rho }\Bigr\}\ \text{and}\
\frac{1}{M}:=\inf_{t\in [1/4,3/4]}\int_{1/4}^{3/4} k(t,s) \,ds.
\end{equation*}
Then $\|Tu\| \geq \|u\|$ for every $u \in K \;\text{with}\; \|u\|=\rho$.
\end{lemma}

\begin{figure}[h]
\centering
\includegraphics[height=5.5cm,angle=0]{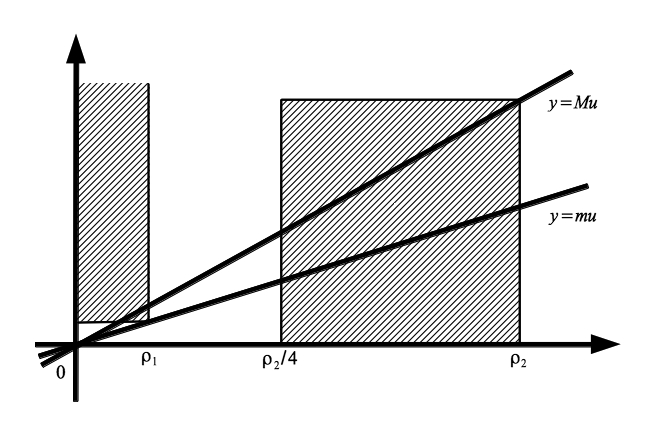}
\caption{The graph of $f$ is not allowed in the shaded region.}
\label{region0}
\end{figure}

\begin{proof}
Take $u\in K$ with $\|u\|=\rho$. For $t\in [1/4,3/4]$ we have
\begin{align*}
Tu(t)=\int_{0}^{1}k(t,s)f(u(s))ds &\geq \int_{1/4}^{3/4}k(t,s)f(u(s))ds\\ &\geq 
\rho f_{{\rho}/{4},{\rho}}\int_{1/4}^{3/4} k(t,s) \,ds\geq \rho=\|u\|.
\end{align*}
\end{proof}
Combining the two Lemmas above we obtain the following Theorem.
\begin{theorem} 
Assume that one of the
following conditions holds.
\begin{itemize}
\item[$(H_{1})$] There exist $\rho_{1},\rho_{2}\in (0, +\infty)$ with
$\rho_{1}<\rho_2/4$ such that $$ f^{0,\rho_{1}}\le
{m}\quad\text{and}\quad f_{\rho_{2}/4,\rho_{2}}\ge {M}.$${}
\item[$(H_{2})$] There exist $\rho_{1},\rho_{2}\in (0, +\infty)$ with
$\rho_{1}<\rho_2$ such that $$ f_{\rho_{1}/4,\rho_{1}}\ge {M}
\quad\text{and}\quad f^{0,\rho_{2}}\le {m}.$${}
\end{itemize}
Then Eq.~\eqref{eqham} has a positive solution in $K$.
\end{theorem}
The case $(H_{1})$ is illustrated in Figure~\ref{region0}.
{}
\begin{example}\label{es1}
Let $\lambda>0$ and consider the BVP
\begin{equation}\label{es1bvp}
u''(t)+\lambda u^{2}(t)=0,\ u(0)=u(1)=0.
\end{equation}
We wish to investigate the values of $\lambda$ for which the BVP~\eqref{es1bvp} admits a non-negative solution of norm less than or equal to $1$.

In this case $m=8$ and $M=16$. By fixing $\rho_2=1$, we have
$$f_{{1}/{4},{1}}=\inf\Bigl\{\frac{f(u)}{1} : 1/4 \leq u \leq {1 }\Bigr\}=f(1/4)=\lambda/16\geq M,$$
if $\lambda\geq 256$. Furthermore, the choice of $\rho_1= 8/\lambda$ gives 
$$
f^{0,{\rho_1}}=\sup\Bigl\{\frac{f(u)}{\rho} : 0 \leq u \leq {\rho_1}\Bigr\}=\frac{f(\rho_1)}{\rho_1}=\lambda \rho_1\leq m.
$$
This implies that, for every $\lambda\geq 256$, the BVP~\eqref{es1bvp} has a non-negative solution $u_{\lambda}$, with $$8/\lambda\leq\|u_{\lambda}\|\leq 1.$$

Note that, by dropping the localization requirement of the solution within the unitary ball, with the same technique it is possible to prove that the BVP~\eqref{es1bvp} admits a non-negative solution for every $\lambda > 0$.

\end{example}
\begin{example}
Consider the BVP
\begin{equation}\label{es2bvp}
u''(t)+\lambda f(u(t))=0,\ u(0)=u(1)=0,
\end{equation}
where $\lambda>0$, $f:[0, +\infty) \to [0, +\infty)$ is continuous, $f_{{1}/{4},{1}}>0$ and $\displaystyle{\lim_{u \to 0^+}}\dfrac{f(u)}{u}=0$. Note that the condition $\lim_{u \to 0^+}\dfrac{f(u)}{u}=0$ implies that  
$
\lambda f^{0,\rho_1} \leq m,
$ for $\rho_1$ sufficiently small.
Then, reasoning as in Example \ref{es1}, it is possible to show that the BVP~\eqref{es2bvp}
 has a positive solution $u_{\lambda}$, with $\|u_{\lambda}\|<1$, for every $\lambda>\dfrac{16}{f_{{1}/{4},{1}}}$. 
\end{example}

\chapter{The fixed point index}\label{fpin}

We now illustrate how to utilize the classical fixed point index in order to prove existence and multiplicity results of solutions for Hammerstein integral equations. The results in this Chapter are essentially based on the manuscripts~\cite{kqljlms, kqljwjde}.

What is the fixed point index of a compact map $T$? Roughly speaking, it is the algebraic count of the fixed points of $T$
in a certain set. The definition is rather technical and
involves the knowledge of the \emph{Leray-Schauder degree}.
Typically the best candidate for a set on which to compute the
fixed point index is a cone.

\begin{proposition}\cite{amann, guolak} Let $D$ be an open bounded set of $X$ with $0 \in D_{K}$ and
$\overline{D}_{K}\ne K$, where $D_{K}=D\cap K$.
Assume that $T:\overline{D}_{K}\to K$ is a compact map such that
$x\neq Tx$ for $x\in \partial D_{K}$. Then the fixed point index
 $i_{K}(T, D_{K})$ has the following properties:
\begin{itemize}
\item[$(1)$] If there exists $e\in K\setminus \{0\}$
such that $x\neq Tx+\lambda e$ for all $x\in \partial D_K$ and all
$\lambda>0$, then $i_{K}(T, D_{K})=0$.
\item[$(1^*)$] If $\|Tx\|\ge \|x\|$ for $x\in
\partial D_K$, then $i_{K}(T, D_{K})=0$.

\item[$(2)$] If $Tx \neq \lambda x$ for all $x\in
\partial D_K$ and all $\lambda > 1$, then $i_{K}(T, D_{K})=1$.

\item[] For example $(2)$ holds if $\|Tx\|\le \|x\|$ for $x\in
\partial D_K$.

\item[(3)] Let $D^{1}$ be open in $X$ such that
$\overline{D^{1}}_{K}\subset D_K$. If $i_{K}(T, D_{K})=1$ and $i_{K}(T,
D_{K}^{1})=0$, then $T$ has a fixed point in $D_{K}\setminus
\overline{D_{K}^{1}}$. The same holds if 
$i_{K}(T, D_{K})=0$ and $i_{K}(T, D_{K}^{1})=1$.
\end{itemize}
\end{proposition}
\begin{definition}
We use the notation
$$K_{\rho}=\{u\in K: \|u\|<\rho\},$$
and we denote by $\partial K_{\rho}$ the boundary relative to $K$.
\end{definition}
\begin{proof}[Proof of Theorem \ref{krasno}]  Assume that $0<r<R$. Then if $T$ has a fixed point on $\partial K_{r}$  or on $ \partial K_{R}$ we are done. Otherwise we have that $i_{K}(T, K_{r})=1$ and $i_{K}(T, K_{R})=0$. By the additivity property of the index we have that $i_{K}(T, K_{R}\setminus \overline{K_{r}})=-1\neq 0$.
Thus there exists a fixed point $u_0$ with $r\leq \|u_0\|\leq R$. The proof of the other case is similar.
\end{proof}
\begin{figure}[h]%\label{fig3}
\centering
\includegraphics[height=5.5cm,angle=0]{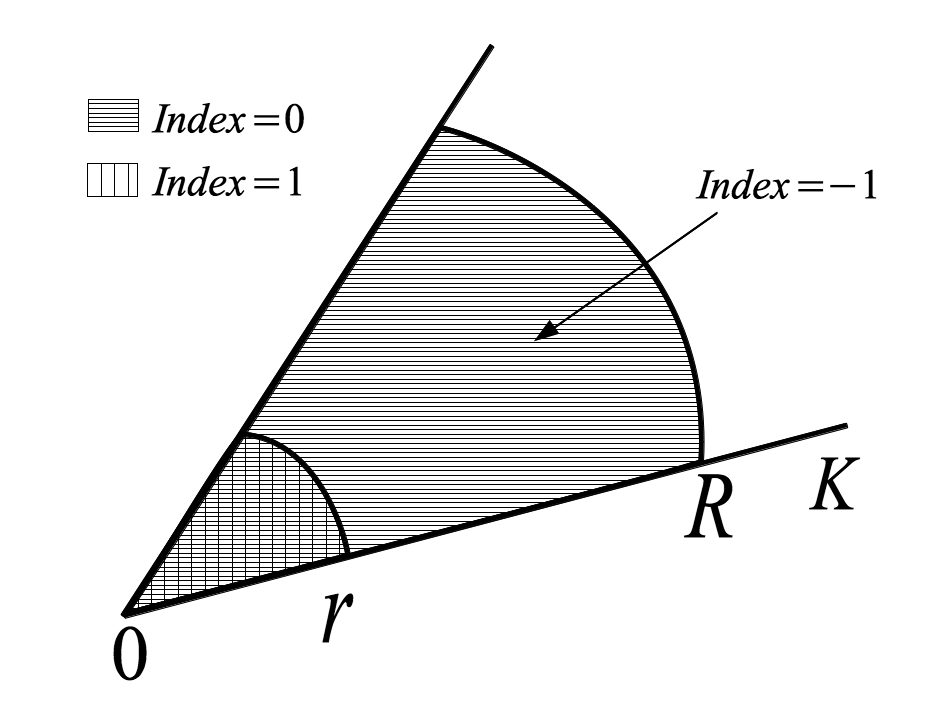}
\caption{A sketch of the proof of Krasnosel'ski\u\i{}'s Theorem.}
\end{figure}{}
We now study Hammerstein integral equations in a slightly more general setting. We assume that the terms that occur in the equation
\begin{equation}\label{eqhamm}
u(t)=\int_{0}^{1}k(t,s)f(u(s))ds :=Tu(t),
\end{equation}
satisfy:
\begin{itemize}
\item  $f:[0, +\infty)\to [0, +\infty)$ is continuous.
\item  $k:[0,1]\times [0,1]\to [0, +\infty)$ is continuous.{}
\item There exist a continuous function  $\Phi:[0,1]\to
[0, +\infty)$, an interval $[a, b]\subset [0,1]$  and a constant $c \in (0,1]$ such that
\begin{align*}
  k(t,s)\leq \Phi(s) \text{ for }& t,s \in [0,1] \text{ and}\\
  c \Phi(s)\leq k(t,s) \text{ for }& t \in [a,b] \text{ and }s \in
  [0,1].
\end{align*}{}
\item $\int_a^b
\Phi(s)\,ds >0$.{}
\end{itemize}

\begin{figure}[h]
\centering
\includegraphics[height=6cm,angle=0]{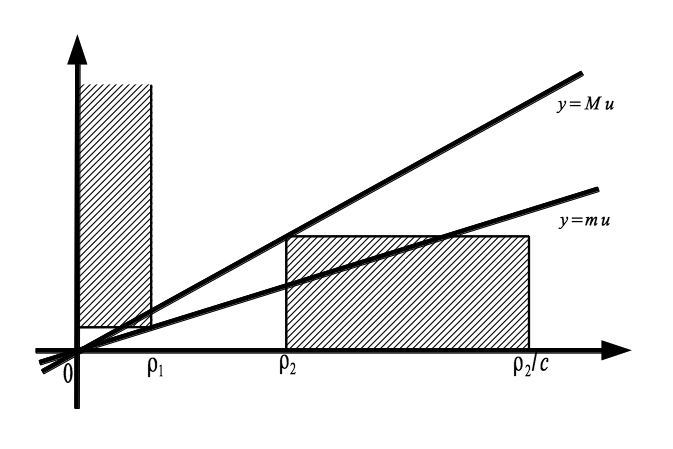}
\caption{One solution.}
\label{region1}
\end{figure}

\begin{figure}[h]
\centering
\includegraphics[height=6cm,angle=0]{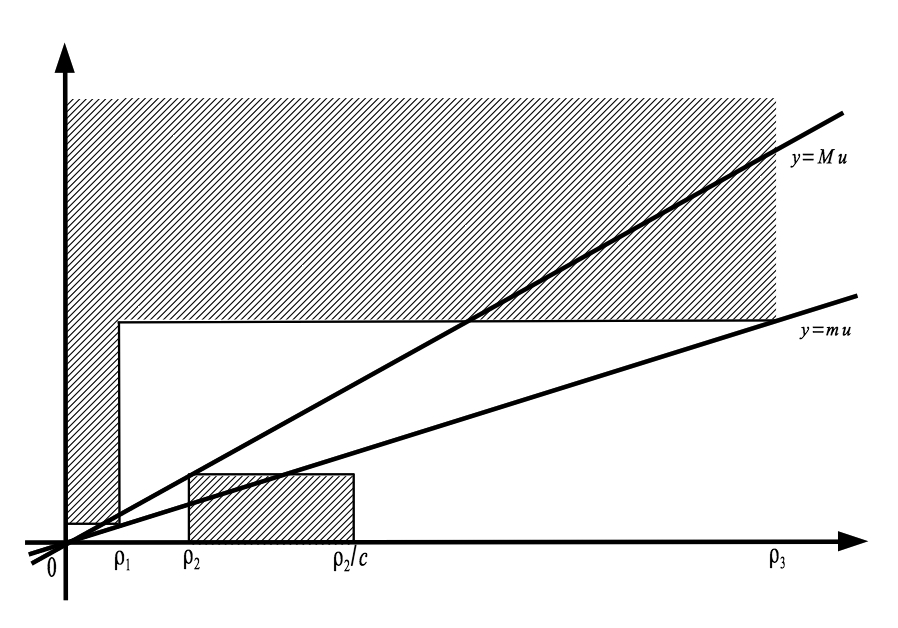}
\caption{Two solutions.}
\label{region2}
\end{figure}

\begin{figure}[h]
\centering
\includegraphics[height=6cm,angle=0]{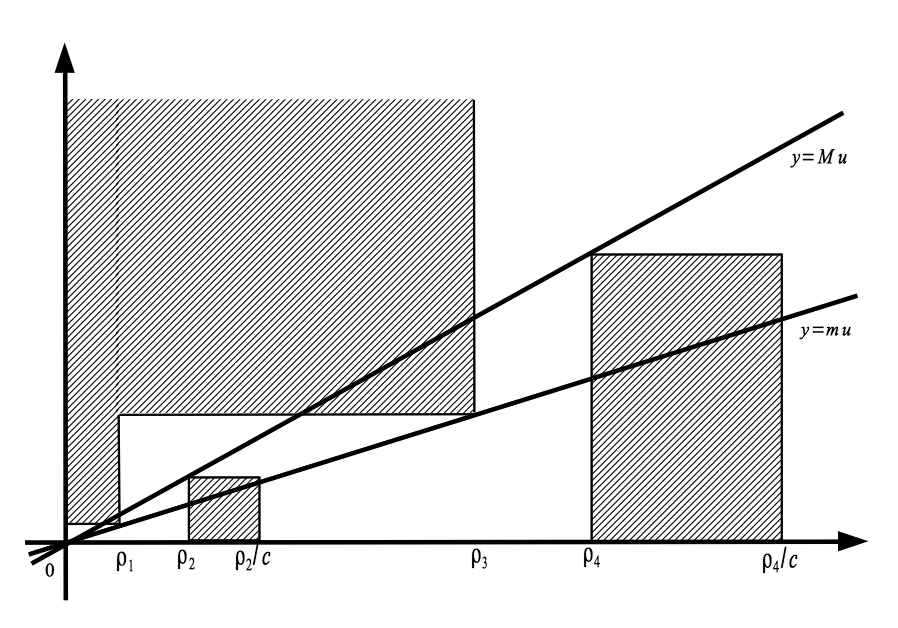}
\caption{Three solutions.}
\label{region3}
\end{figure}

In a similar way as before, we look for fixed points of $T$ in the cone
\begin{equation*}\label{eqcone2}
K=\{u\in C[0,1], u\geq 0: \min_{t \in [a,b]}u(t)\geq c \|u\|\}.
\end{equation*}
It can be shown that, under the hypotheses above, $T$ maps $K$ to $K$ and is compact.
\begin{definition} 
We make use of the set 
\begin{equation*}%\label{vrho}
V_{\rho}=\{u \in K: \displaystyle{\min_{t \in [a,b]}}u(t)<\rho\},
\end{equation*}
The set $V_{\rho}$ was essentially introduced by Lan in~\cite{kqljlms}.
\end{definition}
Note that $K_{\rho}\subset V_{\rho}\subset K_{\rho/c}$.
We now prove two lemmas which give conditions when the fixed point index is either
0 or 1. The line of proof of these results follows, more of less, the one of Lemmas \ref{lemcom} and \ref{lemexp}.
\begin{lemma}%\label{ind1b} 
Assume that 
\begin{enumerate}
\item[$(\mathrm{I}_{\protect\rho }^{1})$]  there exists $\rho> 0$ such that
$
f^{0,\rho} < m,
$
where
\begin{equation*}
  f^{0,{\rho}}:=\sup\Bigl\{\frac{f(u)}{\rho} : 0 \leq u \leq {\rho}\Bigr\}\ \text{and}\
\frac{1}{m}:=\sup_{t\in [0,1]} \int_{0}^{1} k(t,s) \,ds.
\end{equation*}
\end{enumerate}{}
Then the fixed point index, $i_{K}(T,K_{\rho})$, is equal to 1.
\end{lemma}
\begin{proof}
We show that $\mu u \neq Tu$ for every $u \in \partial K_{\rho }$
and for every $\mu \geq 1$. 
In fact, if this does not happen, there exist $\mu \geq 1$ and $u\in
\partial K_{\rho }$ such that $\mu u=Tu$,
that is
$$\mu u(t)= \int_{0}^{1}
k(t,s)f(u(s))\,ds,$$
Taking  the supremum for $t\in [0,1]$ gives
\begin{equation*}
\mu \rho\leq 
\sup_{t\in [0,1]}\int_{0}^{1}k(t,s)f(u(s))\,ds
 \leq\rho f^{0,\rho}\cdot\sup_{t\in [0,1]}\int_{0}^{1}k(t,s)\,ds
 <\rho.
\end{equation*}

This
contradicts the fact that $\mu \geq 1$ and proves the result.
\end{proof}
\begin{lemma}
\label{idx0b1} Assume that

\begin{enumerate}
\item[$(\mathrm{I}_{\protect\rho }^{0})$] there exist $\rho >0$ such that
such that
$
{f_{\rho ,{\rho /c}}}>M(a,b),
$
where
$$ 
f_{\rho ,{\rho /c}} =\inf \left\{\frac{f(u)}{\rho }%
:\;\rho\leq u \leq \rho /c\right\}
\ \text{and}\ \frac{1}{M(a,b)} :=\inf_{t\in [a,b]}\int_{a}^{b}k(t,s)\,ds.$$
\end{enumerate}

Then $i_{K}(T,V_{\rho})=0$.
\end{lemma}
\begin{proof} Let $e(t)\equiv 1$, then $e \in K$. 
We prove that
\begin{equation*}
u\ne Tu+\lambda e\quad\text{for  all } u\in \partial
V_{\rho}\text{ and } \lambda \geq 0.
\end{equation*}
In fact, if not, there exist $u\in \partial V_\rho$ and $\lambda\geq 0$ such that $u=Tu+\lambda e$. Then we have$$
u(t)=\int_{0}^{1}
k(t,s)f(u(s))\,ds+\lambda.$$
Thus we get, for $t\in[a,b]$,
\begin{equation*}
u(t)=\int_{0}^{1}
k(t,s)f(u(s))\,ds+ \lambda \ge
\int_{a}^{b}
k(t,s)f(u(s))\,ds\ge
\rho f_{\rho ,{\rho /c}}
\;  
\int_{a}^{b}
k(t,s)\,ds.
\end{equation*}
Taking the minimum over $[a,b]$ gives
$\rho>\rho$ a contradiction.
\end{proof}
\begin{remark}\label{comp}
In order to compare the two approaches, proving the index zero result by means of the condition $\|Tu\| \geq \|u\|$ for every $u \in \partial K_{\rho/c}$, as in the application of the Krasnosel'ski\u\i{} Theorem,  would require
$$f(u)\geq M\rho/c,\ \text{for every}\ u\in  [\rho ,\rho /c],$$ a more stringent requirement. 
\end{remark}
\begin{remark}
Note also that we used strict inequalities in the conditions  $(\mathrm{I}_{\rho}^{0})$ and  $(\mathrm{I}_{\rho}^{1})$. This fact is particularly convenient for proving the existence of multiple solutions, this is done in the following Theorem.
\end{remark}
\begin{theorem}
\label{thmmsol1} The integral equation~\eqref{eqhamm} has at least one positive solution
in $K$ if either of the following conditions holds.

\begin{enumerate}

\item[$(S_{1})$] There exist $\rho _{1},\rho _{2}\in (0, +\infty )$ with $\rho
_{1}/c<\rho _{2}$ such that $(\mathrm{I}_{\rho _{1}}^{0})$ and $(\mathrm{I}_{\rho _{2}}^{1})$ hold.

\item[$(S_{2})$] There exist $\rho _{1},\rho _{2}\in (0, +\infty )$ with $\rho
_{1}<\rho _{2}$ such that $(\mathrm{I}_{\rho _{1}}^{1})$ and $(\mathrm{I}%
_{\rho _{2}}^{0})$ hold.
\end{enumerate}
The integral equation \eqref{eqhamm} has at least two positive solutions in $K$ if one of
the following conditions holds.

\begin{enumerate}

\item[$(S_{3})$] There exist $\rho _{1},\rho _{2},\rho _{3}\in (0, +\infty )$
with $\rho _{1}/c<\rho _{2}<\rho _{3}$ such that $(\mathrm{I}_{\rho
_{1}}^{0}),$ $(
\mathrm{I}_{\rho _{2}}^{1})$ $\text{and}\;\;(\mathrm{I}_{\rho _{3}}^{0})$
hold.

\item[$(S_{4})$] There exist $\rho _{1},\rho _{2},\rho _{3}\in (0, +\infty )$
with $\rho _{1}<\rho _{2}$ and $\rho _{2}/c<\rho _{3}$ such that $(\mathrm{I}%
_{\rho _{1}}^{1}),\;\;(\mathrm{I}_{\rho _{2}}^{0})$ $\text{and}\;\;(\mathrm{I%
}_{\rho _{3}}^{1})$ hold.
\end{enumerate}
The integral equation \eqref{eqhamm} has at least three positive solutions in $K$ if one
of the following conditions holds.

\begin{enumerate}
\item[$(S_{5})$] There exist $\rho _{1},\rho _{2},\rho _{3},\rho _{4}\in
(0, +\infty )$ with $\rho _{1}/c<\rho _{2}<\rho _{3}$ and $\rho _{3}/c<\rho
_{4}$ such that $(\mathrm{I}_{\rho _{1}}^{0}),$ $(\mathrm{I}_{\rho _{2}}^{1}),\;\;(\mathrm{I}%
_{\rho _{3}}^{0})\;\;\text{and}\;\;(\mathrm{I}_{\rho _{4}}^{1})$ hold.

\item[$(S_{6})$] There exist $\rho _{1},\rho _{2},\rho _{3},\rho _{4}\in
(0, +\infty )$ with $\rho _{1}<\rho _{2}$ and $\rho _{2}/c<\rho _{3}<\rho _{4}$
such that $(\mathrm{I}_{\rho _{1}}^{1}),\;\;(\mathrm{I}_{\rho
_{2}}^{0}),\;\;(\mathrm{I}_{\rho _{3}}^{1})$ $\text{and}\;\;(\mathrm{I}%
_{\rho _{4}}^{0})$ hold.
\end{enumerate}
\end{theorem}
\begin{proof}
We sketch the proof in Figures~\ref{onesolA},  \ref{onesolB}, \ref{twosolA}, and \ref{twosolB}.

\begin{figure}[ht]
\begin{minipage}[b]{0.45\linewidth}
\centering
\includegraphics[height=5.1cm,angle=0]{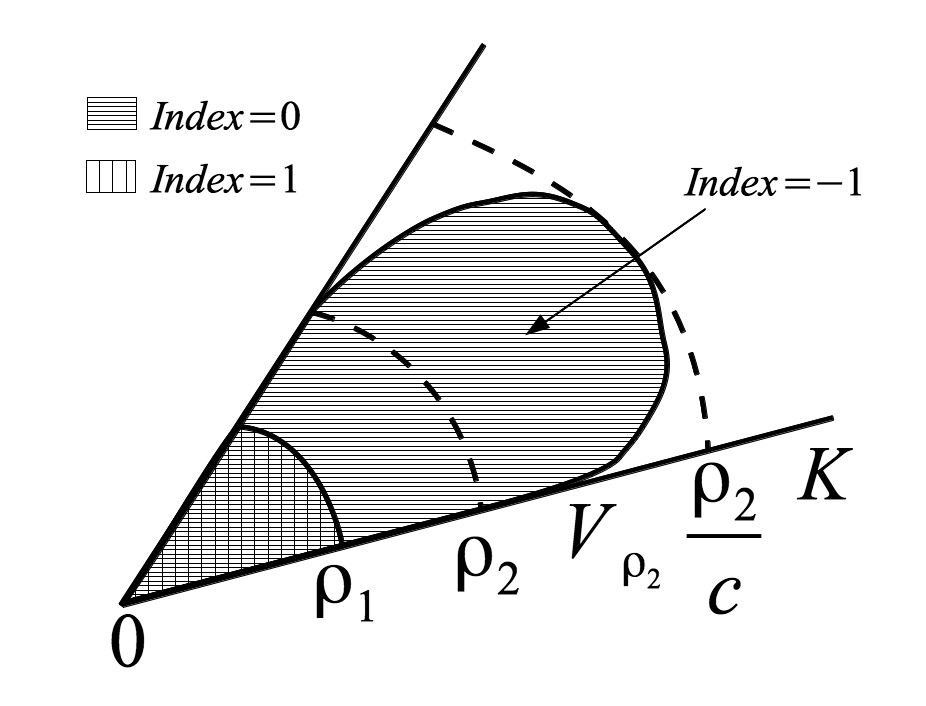}
\caption{One positive solution.}
\label{onesolA}
\end{minipage}
\hspace{0.5cm}
\begin{minipage}[b]{0.45\linewidth}
\centering
\includegraphics[height=5.1cm,angle=0]{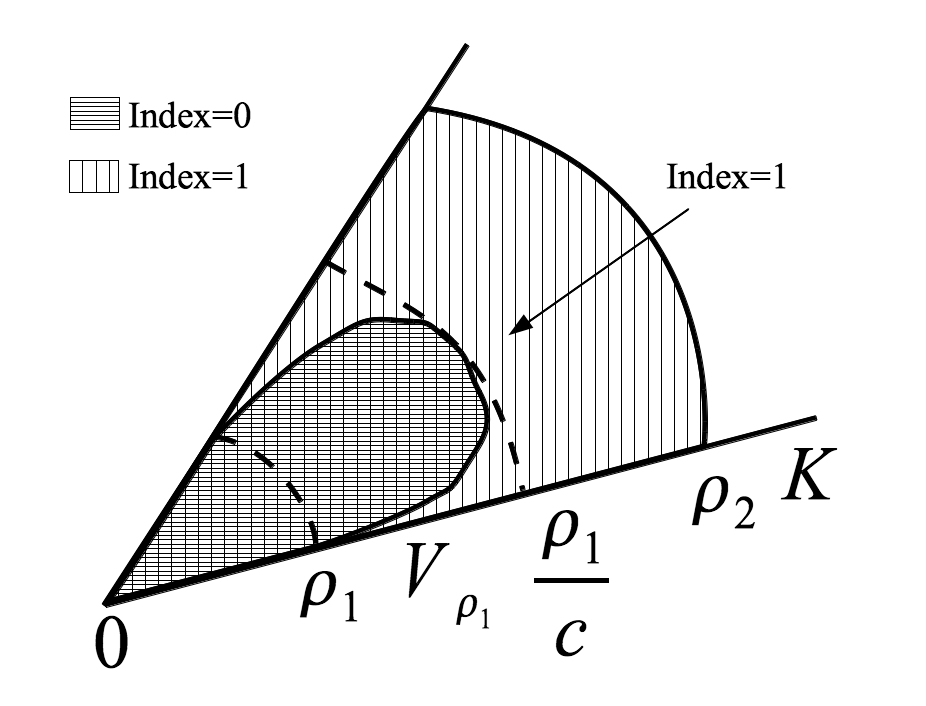}
\caption{One positive solution.}
\label{onesolB}
\end{minipage}
\end{figure}

\begin{figure}[ht]
\begin{minipage}[b]{0.45\linewidth}
\centering
\includegraphics[height=5.1cm,angle=0]{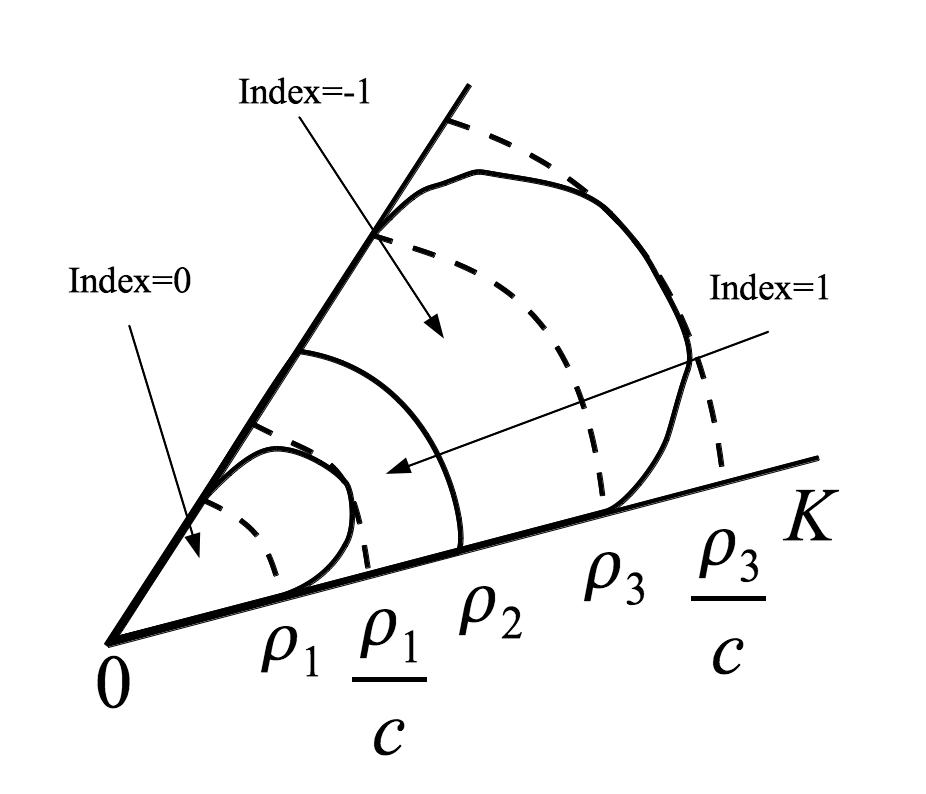}
\caption{Two positive solutions.}
\label{twosolA}
\end{minipage}
\hspace{0.5cm}
\begin{minipage}[b]{0.45\linewidth}
\centering
\includegraphics[height=5.1cm,angle=0]{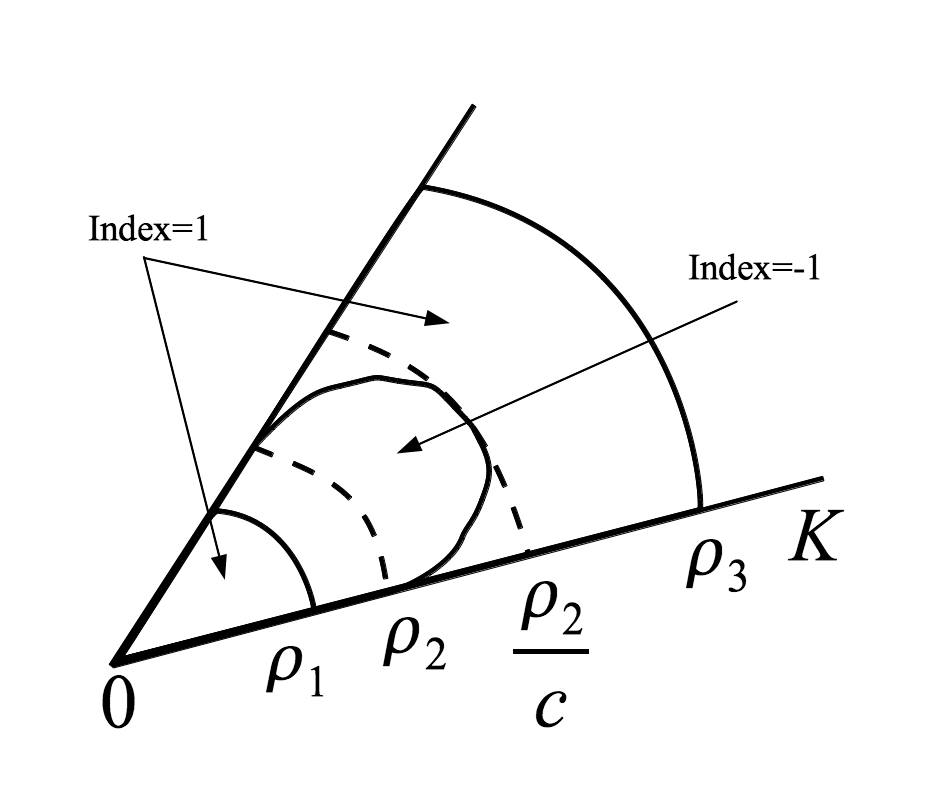}
\caption{Two positive solutions.}
\label{twosolB}
\end{minipage}
\end{figure}
\end{proof}
\begin{example}
In view of Remark~\ref{comp} is clear that we improve the growth assumptions on the nonlinearity occurring in the BVP~\eqref{EQ-dir}-\eqref{BC-dir} since we need 
$$f(u)\geq 16 \rho,\ \text{for every}\ u\in  [\rho ,\rho /c],$$ instead of $$f(u)\geq 64\rho,\ \text{for every}\ u\in  [\rho ,\rho /c].$$
\end{example}
\begin{example}
We can apply Theorem~\ref{thmmsol1} to study the existence of positive solutions for the BVP
\begin{align}\label{BVP2}
u''(t)+f(u(t))=0,\ u(0)=u'(1)=0.
\end{align}
The BCs in~\eqref{BVP2} are called~\emph{right focal} BCs or, sometimes, \emph{mixed} BCs, since are on the left side of the interval of Dirichlet type and on the other side of Neumann type.

In order to construct the Green's function we consider the linear problem
\begin{align*}%\label{linBVP2}
u''(t)+y(t)=0,\ u(0)=u'(1)=0.
\end{align*}
If we integrate $u''+y=0$ we obtain
$$u'(t)=u'(0)+\int_{0}^{t}-y(s)ds$$ 
and, using the BC $u'(1)=0$ we get
$$u'(0)=\int_{0}^{1}y(s)ds$$ 
and
$$u(t)=u(0)+tu'(0)-\int_{0}^{t}\int_{0}^{w}y(s)dsdw.$$
By using the BCs and the Cauchy formula for iterated kernels, we obtain 
$$u(t)=t\int_{0}^{1}y(s)ds-\int_{0}^{t}(t-s)y(s)ds.$$
This gives
$$
u(t)=\int_{0}^{1}k(t,s)y(s)ds,
$$
where
\begin{equation}\label{kRob}
k(t,s)=\begin{cases}
  s & \text{if}\;s \leq t,\\
  t &\text{if}\; s>t.
\end{cases}
\end{equation}
Therefore the solution of the BVP~\eqref{BVP2} is  given by 
\begin{equation*}%\label{eqhamm2}
u(t)=\int_{0}^{1}k(t,s)f(u(s))ds :=Tu(t),
\end{equation*}
In this case one may take as an upper bound for the kernel $\Phi(s)=s$ and show that $k(t,s)\geq as$ on $[a,b]\times [0,1]$. Thus $[a,b]$ can be chosen arbitrarily in $(0,1]$. In this case  $m=2$ and the choice of $[a,b]=[1/2,1]$ gives $c=1/2$ and $M=4$. Note that this choice for $[a,b]$ is optimal in the sense that provides the minimal $M$ to be satisfied in condition $(\mathrm{I}_{\rho }^{0})$.

In Figures~\ref{region1}, \ref{region2} and \ref{region3} we illustrate the allowed growth of a nonlinearity $f$ for the existence of one, two and three positive solutions, that correspond to the cases $(S_{2})$, $(S_{4})$, $(S_{6})$ of Theorem~\ref{thmmsol1}.
\end{example}

\section{A non-existence result}
We now prove a simple non-existence result for the integral equation~\eqref{eqhamm}.
\begin{theorem}\label{nonext}
 Assume that one of the following conditions holds: 
 \begin{enumerate}
\item[$(1)$] $f(u)<m u$ for $u >0$,
\item[$(2)$] $f(u)>M u$ for $u >0$.
\end{enumerate}
Then the equation~\eqref{eqhamm} has no non-trivial solution in $K$. 
\end{theorem}

\begin{proof}
$(1)$ Assume, on the contrary, that there exists $u\in K$, $u\not\equiv 0$ such that $u=Tu$ and let $t_0\in [0,1]$ be such that $\|u\|=u(t_0)$. Then we have
\begin{align*}
\|u\|=u(t_0)=\int_{0}^{1}k(t_0,s)f(u(s))\,ds <& \int_{0}^{1}k(t_0,s) m u(s)\, ds\\
 \leq &  m\|u\| \Bigl( \int_{0}^{1}k(t_0,s)\, ds\Bigr)\leq \|u\|
\end{align*}
a contradiction.\par
$(2)$ Assume, on the contrary, that there exists $u\in K$, $u\not\equiv 0$ such that $u=Tu$ and let $\eta\in [a,b]$ be such that $u(\eta)=\min_{t\in[a,b]}u(t)$. For $t\in [a,b]$ we have
$$
u(t)=  \int_0^1 k(t,s)f(u(s))\,ds\geq \int_{a}^{b} k(t,s)f(u(s))\,ds>  M \int_{a}^{b} k(t,s)u(s)\,ds.
$$
Taking the infimum  for $t\in [a,b]$, we have
$$
\min_{t\in[a,b]}u(t)> M \inf_{t\in[a,b]}\int_{a}^{b} k(t,s)u(s)\,ds.
$$
Thus we obtain
$$
u(\eta)>M u(\eta) \inf_{t\in[a,b]}\int_{a}^{b} k(t,s)\,ds \geq u(\eta),
$$
a contradiction.
\end{proof}
\begin{remark}
Note that since $T$ maps $P$ into $K$, we have proven that~\eqref{eqhamm} has no non-trivial solution in $P$. In Figures \ref{region-no1} and \ref{region-no2} we illustrate the growths of a nonlinearity $f$ in the case of Theorem~\ref{nonext}.
\begin{figure}[ht]
\centering
\includegraphics[height=5.5cm,angle=0]{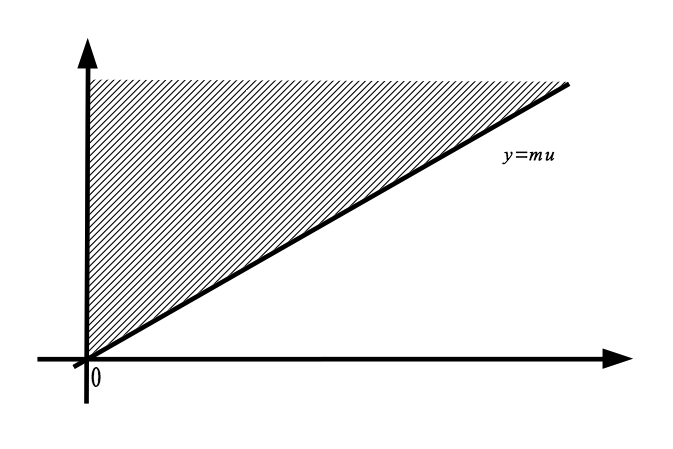}
\caption{Non-existence in $K$, $f$ ``small''.}
\label{region-no1}
\end{figure}
\begin{figure}[ht]
\centering
\includegraphics[height=5.5cm,angle=0]{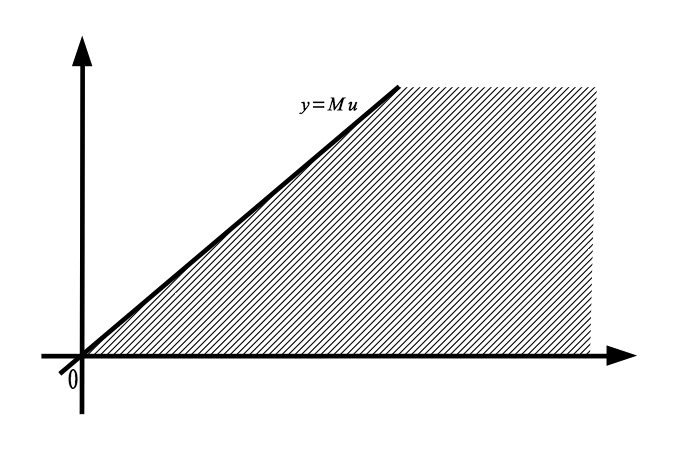}
\caption{Non-existence in $K$, $f$ ``large''.}
\label{region-no2}
\end{figure}
\end{remark}
\chapter{Nonnegative solutions of systems of BVPs}\label{secsyst}
We now  discuss the existence of non-negative solutions for the system of second order BVPs
\begin{gather}
\begin{aligned} \label{sys}
&u^{\prime \prime }(t)+f_{1}(u(t),v(t))=0,\ t \in (0,1), \\
&v^{\prime \prime }(t)+f_{2}(u(t),v(t))=0,\ t \in (0,1), \\
&u(0)= u(1)=v(0)=v'(1)=0.%
\end{aligned}%
\end{gather}
The results in this Chapter are essentially based on the manuscripts~\cite{df-gi-do, gipp-ns, gipp-nonlin}.

In similar manner as the case of one equation, we would like to use a formulation that involves integral equations.
In particular we rewrite the system~\eqref{sys} as a system of Hammerstein integral equations, that is
\begin{gather}
\begin{aligned}\label{ham-sys}
u(t)=\int_{0}^{1}k_1(t,s)f_1(u(s),v(s))\,ds, \\
v(t)= \int_{0}^{1}k_2(t,s)f_2(u(s),v(s))\,ds,
\end{aligned}
\end{gather}
where $k_1$ is given by~\eqref{kDir} and $k_2$ is given by~\eqref{kRob}.

We assume the following.
\begin{itemize}
\item For every $i=1,2$, $f_i:  [0, +\infty)\times [0, +\infty) \to [0, +\infty)$ is continuous.
{}
\end{itemize}

We work in  the space $C[0,1] \times C[0,1]$ endowed (with abuse of notation) with the norm
$$\|(u,v)\|:=\max\{\|u\|, \|v\| \}.$$ 
Let
$$
\tilde{K_i}:= \{ w \in C[0,1]: w(t)\geq 0\ \text{and}\; \min_{t\in [a_i,b_i]} w(t) \geq c_i
\|w\|  \},
$$
where $[a_1,b_1]=[1/4,3/4]$, $[a_2,b_2]=[1/2,1]$, $c_1=1/4$ and $c_2=1/2$,
 and consider the
cone $K$ in $C[0,1] \times C[0,1]$ defined by
\begin{equation*}
\label{cone}
\begin{array}{c}
 K:=  \{ (u,v) \in \tilde{K_1} \times \tilde{K_2}  \}.
\end{array}
\end{equation*}

For a \emph{positive} solution of the system  \eqref{sys}
we mean a solution $(u,v)\in K$ of \eqref{ham-sys} such
that $\|(u,v)\|>0$.

Under our assumptions, a routine check shows that the integral
operator
\begin{equation*}
T \Bigl(
\begin{array}{c}
u\\
v
\end{array}
\Bigr)(t)
:=
\left(
\begin{array}{c}
T_1(u,v)(t)\\
T_2(u,v)(t)
\end{array}
\right),
\end{equation*}
where
$$
T_i(u,v)(t):=\int_{0}^{1}k_i(t,s)f_i(u(s),v(s))\,ds,
$$
leaves $K$ invariant and is compact.

For our fixed point index calculations we work with the following (relative) open bounded sets in $K$:
$$
K_{\rho} =  \{ (u,v) \in K : \|(u,v)\|< \rho  \},
$$
and
$$
V_\rho=\{(u,v) \in K: \min_{t\in [a_1,b_1]}u(t)<\rho\ \text{and}\   \min_{t\in [a_2,b_2]}v(t)<\rho\}.
$$
Set $c=\min\{{c_1},{c_2}\}$. The set $V_\rho$ (in the context of systems) was introduced in~\cite{gipp-ns} and is equal to the set
called $\Omega^{\rho /c}$ in~\cite{df-gi-do}. $\Omega^{\rho /c}$ is an extension to the case of systems of a set given by Lan~\cite{kqljlms}. 
As before
we denote by $\partial K_{\rho}$ and $\partial V_{\rho}$
 the boundary of $K_{\rho}$ and $V_{\rho}$ relative to $K$.
 
The following Lemma provides some useful properties of the set $V_{\rho}$.

\begin{lemma}  %\label{esca} 
  The sets defined above have the following properties:
\begin{itemize}
\item $K_{\rho}\subset V_{\rho}\subset K_{\rho/c}$.
\item $(w_1,w_2) \in \partial V_{\rho}$ \; iff \; $(w_1,w_2)\in K$ and $\displaystyle
\min_{t\in [a_i,b_i]} w_i(t)= \rho$ for some $i\in \{1,2\}$ and $\displaystyle
\min_{t\in [a_j,b_j]}w_j(t)
\le \rho$ for $j\neq i$.
\item If $(w_1,w_2) \in \partial V_{\rho}$, then for some $i\in\{1,2\}$ $\rho \le w_i(t) \le \rho/c$
 for each $t \in [a_i,b_i]$ and for $j\neq i$ we have $0 \leq w_j(t) \leq \rho/c$ for each $t\in [a_j,b_j]$ and $\|w_j\| \leq \rho/c$.
\end{itemize}
\end{lemma}

We can now provide some index results for the case of systems. 
 
\begin{lemma} \label{ind1b-sys}
Assume that
\begin{enumerate}
\item[$(\mathrm{I}^{1}_{\rho})$]%\label{EqB} 
there exists $\rho>0$ such that, for every $i=1,2$,
$
 f_i^{0,\rho}<m_i,
$
where
\begin{equation*}
f_i^{0,{\rho}}=\sup \Bigl\{ \frac{f_i(u,v)}{ \rho}:\; (u,v)\in
[0, \rho]\times[0, \rho]\Bigr\}\ \text{and}\
\frac{1}{m_i}=\sup_{t\in [0,1]} \int_{0}^{1} k_i(t,s) \,ds.
\end{equation*}{}
\end{enumerate}
Then $i_{K}(T,K_{\rho})$ is equal to 1.
\end{lemma}

\begin{proof}
We show that $\mu (u,v)\neq T(u,v)$ for every $(u,v) \in \partial
K_{\rho}$ and for every $\mu \geq 1$;
this ensures that the index is 1 on $K_{\rho}$. In fact, if this does not happen, there exists
$\mu \geq 1$ and $(u,v) \in
\partial K_{\rho}$ such that $\mu (u,v)= T(u,v)$.
Assume, without loss of generality, that $\|u\|=\rho$ and $\|v\|\leq\rho$.
Then
$$
 \mu
u(t)= T_1(u,v)(t)=\int_{0}^{1}k_1(t,s)f_1(u(s),v(s))\,ds.
$$
Taking  the supremum for $t\in [0,1]$ gives
\begin{equation*}
\mu \rho\leq 
\sup_{t\in [0,1]}\int_{0}^{1}k_1(t,s)f_1(u(s),v(s))\,ds
 \leq\rho f_1^{0,{\rho}}\cdot\sup_{t\in [0,1]}\int_{0}^{1}k_1(t,s)\,ds
 <\rho.
\end{equation*}
This
contradicts the fact that $\mu \geq 1$ and proves the result.
\end{proof}
We give a first Lemma that shows that the index is 0 on a set $V_{\rho}$.  
\begin{lemma}\label{idx0b1-sys}
Assume that
\begin{enumerate}
\item[$(\mathrm{I}^{0}_{\rho})$]  there exist $\rho>0$ such that, for every $i=1,2$,
$
f_{i,(\rho, \rho/c)}>M_i,
$
where
\begin{multline*}
f_{1,(\rho,{\rho / c})}=\inf \Bigl\{ \frac{f_1(u,v)}{ \rho}:\; (u,v)\in [\rho,\rho/c]\times[0, \rho/c]\Bigr\},\\
f_{2,(\rho,{\rho / c})}=\inf \Bigl\{ \frac{f_2(u,v)}{ \rho}:\; (u,v)\in [0,\rho/c]\times[\rho, \rho/c]\Bigr\}\\
 \text{and}\ \frac{1}{M_i}=\inf_{t\in
[a_i,b_i]}\int_{a_i}^{b_i} k_i(t,s) \,ds.
\end{multline*}
\end{enumerate}
Then $i_{K}(T,V_{\rho})=0$.
\end{lemma}

\begin{proof}
Let $e(t)\equiv 1$ for $t\in [0,1]$. Then $(e,e)\in K$. We prove that
\begin{equation*}
(u,v)\ne T(u,v)+\mu (e,e)\quad\text{for } (u,v)\in \partial
V_{\rho}\quad\text{and } \mu \geq 0.
\end{equation*}
In fact, if this does not happen, there exist $(u,v)\in \partial V_\rho$ and
$\mu \geq 0$ such that $(u,v)=T(u,v)+\mu (e,e)$.
Without loss of generality, we can assume that for all $t\in [a_1,b_1]$ we have
$$
\rho\leq u(t)\leq {\rho/c},\\\ \min u(t)=\rho \\\ \text{and  }\\\ 0\leq v(t)\leq {\rho/c}.
$$
Then, for $t\in [a_1,b_1]$, we obtain
\begin{multline*}
u(t)= \int_{0}^{1} k_1(t,s)f_1(u(s),v(s))ds
+ \mu e\\
\geq \int_{a_1}^{b_1}k_1(t,s)f_1(u(s),v(s))\,ds+{\mu}\ge
\rho f_{1,(\rho,{\rho / c})}
\;  
\int_{a_1}^{b_1}k_1(t,s)\,ds +{\mu}.
\end{multline*}
Thus, we obtain
$\rho=\min_{t\in [a_1,b_1]}u(t)>
\rho+\mu\geq \rho$,
a contradiction.
\end{proof}

The following Lemma shows that the index is 0 on
$V_{\rho}$; this time we have to control the growth of just one
nonlinearity $f_i$, at the cost of having to deal with a larger domain. This allows to deal with nonlinearities with different
growth, see also the papers~\cite{gifmpp-cnsns, gipp-nonlin,
precup1, precup2, ya1}.
\begin{lemma}\label{idx0b3-sys}
Assume that
\begin{enumerate}
\item[$(\mathrm{I}^{0}_{\rho})^{\star}$] there exist $\rho>0$ such that, for some $i=1,2$,
$
f^*_{i,(0, \rho/c)} >M_i.$
\end{enumerate}
where
\begin{equation*}
f^*_{i,(0,{\rho / c})}=\inf \Bigl\{ \frac{f_i(u,v)}{ \rho}:\; (u,v)\in [0,\rho/c]\times[0, \rho/c]\Bigr\}.
\end{equation*}
Then $i_{K}(T,V_{\rho})=0$.
\end{lemma}
\begin{proof}
Suppose that the condition $(\mathrm{I}^{0}_{\rho})^{\star}$ holds for $i=1$.
Let $e(t)\equiv 1$ for $t\in [0,1]$. Then $(e,e)\in K$. We prove that
\begin{equation*}
(u,v)\ne T(u,v)+\mu (e,e)\quad\text{for } (u,v)\in \partial
V_{\rho}\quad\text{and } \mu \geq 0.
\end{equation*}
In fact, if this does not happen, there exist $(u,v)\in \partial V_\rho$ and
$\mu \geq 0$ such that $(u,v)=T(u,v)+\mu (e,e)$.
So, for all $t\in [a_1,b_1]$, $\min u(t)\leq \rho$ and for $t\in [a_2,b_2]$, $\min v(t)\leq \rho$.
We have, for $t\in [0,1]$,
$$
u(t)=\int_{0}^{1} k_1(t,s)f_1(u(s),v(s))ds
+ \mu e
$$
and, reasoning as in the proof of Lemma \ref {idx0b1-sys}, we obtain  $\rho\geq \min_{t\in [a_1,b_1]}u(t) >\rho+\mu\geq \rho,$ 
a contradiction.
\end{proof}
\begin{theorem}\label{mult-sys}
The system \eqref{ham-sys} has at least one positive solution
in $K$ if either of the following conditions holds.

\begin{enumerate}

\item[$(S_{1})$] There exist $\rho _{1},\rho _{2}\in (0, +\infty )$ with $\rho
_{1}/c<\rho _{2}$ such that $(\mathrm{I}_{\rho _{1}}^{0})\;\;[\text{or}\;(%
\mathrm{I}_{\rho _{1}}^{0})^{\star }],$ $(\mathrm{I}_{\rho _{2}}^{1})$ hold.

\item[$(S_{2})$] There exist $\rho _{1},\rho _{2}\in (0, +\infty )$ with $\rho
_{1}<\rho _{2}$ such that $(\mathrm{I}_{\rho _{1}}^{1}),\;\;(\mathrm{I}%
_{\rho _{2}}^{0})$ hold.
\end{enumerate}

The system \eqref{ham-sys} has at least two positive solutions in $K$ if one of
the following conditions holds.

\begin{enumerate}

\item[$(S_{3})$] There exist $\rho _{1},\rho _{2},\rho _{3}\in (0, +\infty )$
with $\rho _{1}/c<\rho _{2}<\rho _{3}$ such that $(\mathrm{I}_{\rho
_{1}}^{0})$ $[\text{or}\;(\mathrm{I}_{\rho _{1}}^{0})^{\star }],$ $(
\mathrm{I}_{\rho _{2}}^{1})$ $\text{and}\;\;(\mathrm{I}_{\rho _{3}}^{0})$
hold.

\item[$(S_{4})$] There exist $\rho _{1},\rho _{2},\rho _{3}\in (0, +\infty )$
with $\rho _{1}<\rho _{2}$ and $\rho _{2}/c<\rho _{3}$ such that $(\mathrm{I}%
_{\rho _{1}}^{1}),\;\;(\mathrm{I}_{\rho _{2}}^{0})$ $\text{and}\;\;(\mathrm{I%
}_{\rho _{3}}^{1})$ hold.
\end{enumerate}

The system \eqref{ham-sys} has at least three positive solutions in $K$ if one
of the following conditions holds.

\begin{enumerate}
\item[$(S_{5})$] There exist $\rho _{1},\rho _{2},\rho _{3},\rho _{4}\in
(0, +\infty )$ with $\rho _{1}/c<\rho _{2}<\rho _{3}$ and $\rho _{3}/c<\rho
_{4}$ such that $(\mathrm{I}_{\rho _{1}}^{0})\;\;[\text{or}\;(\mathrm{I}%
_{\rho _{1}}^{0})^{\star }],$ $(\mathrm{I}_{\rho _{2}}^{1}),\;\;(\mathrm{I}%
_{\rho _{3}}^{0})\;\;\text{and}\;\;(\mathrm{I}_{\rho _{4}}^{1})$ hold.

\item[$(S_{6})$] There exist $\rho _{1},\rho _{2},\rho _{3},\rho _{4}\in
(0, +\infty )$ with $\rho _{1}<\rho _{2}$ and $\rho _{2}/c<\rho _{3}<\rho _{4}$
such that $(\mathrm{I}_{\rho _{1}}^{1}),\;\;(\mathrm{I}_{\rho
_{2}}^{0}),\;\;(\mathrm{I}_{\rho _{3}}^{1})$ $\text{and}\;\;(\mathrm{I}%
_{\rho _{4}}^{0})$ hold.
\end{enumerate}
\end{theorem}
The proof follows as the one of Theorem~\ref{thmmsol1} and is omitted.
\begin{remark}
In Lemmas \ref{ind1b-sys}, \ref{idx0b1-sys}, \ref{idx0b3-sys} and in Theorem~\ref{mult-sys} we used, for simplicity, the same radii for the component $u$ and $v$. The reader might find different radii in the components in the manuscripts~\cite{chzh2, gipp-nodea}. 
A non-existence result, similar to Theorem~\ref{nonext} can be stated in the case of systems, we refer the reader to~\cite{gipp-nodea}. 
\end{remark}
\begin{example}
Consider the BVP
\begin{gather}
\begin{aligned} \label{sys-ex}
&u^{\prime \prime }(t)+18+\sin(u(t)v(t))=0,\ t \in (0,1), \\
&v^{\prime \prime }(t)+e^{\frac{u^2(t)+v^2(t)}{25}}-1=0,\ t \in (0,1), \\
&u(0)= u(1)=v(0)=v'(1)=0.%
\end{aligned}%
\end{gather}
We  show that 
 $f_1$ satisfies conditions $(\mathrm{I}^{0}_{1})^{\star}$ and $(\mathrm{I}^{1}_{5})$, while 
$f_2$ satisfies $(\mathrm{I}^{1}_{5})$. We take $[a_1,b_1]=[1/4,3/4]$ and $[a_2,b_2]=[1/2,1]$, thus $c_{1}=1/4$, $c_{2}=1/2$ and 
$c=\min\{1/4,1/2\}=1/4$. In this case $m_1=8, \, M_1=16, \,m_2=2, \, M_2=4$, and for $(u,v)\in [0,5]\times[0,5]$ we have
\begin{equation*}
\begin{array}{c}
f_{1}(u,v)=18+\sin(uv)\leq 19 < 8 \times 5=40,\\
f_{2}(u,v)=e^{\frac{u^2+v^2}{25}}-1\leq e^{2}-1  < 2 \times 5=10,
\end{array}
\end{equation*}
so that condition $(\mathrm{I}^{1}_{5})$ holds.

Furthermore, for $(u,v)\in [0,4]\times[0,4]$, we have
$$
f_{1}(u,v)=18+\sin(uv)\geq 17> 16 \cdot 1, 
$$
so that condition  $(\mathrm{I}^{0}_{1})^{\star}$ holds. Thus condition $(S_1)$ of Theorem \ref{mult-sys} is satisfied, providing the existence of at least 
one positive solution $(u_0,v_0)$ of the system \eqref{sys-ex} and, furthermore, we have that 
$1\leq \|(u_0,v_0)\|\leq 4.$
\end{example}

\chapter{More general BCs}\label{mgenbcs}
We now move to the case of non-homogeneous BCs and illustrate how the machinery developed in the previous Chapters can be adapted to this new setting.
\section{A three-point problem}
We begin with a simple three-point problem, by considering the ODE 
\begin{equation}\label{eqD}
u''(t)+f(u(t))=0,\ \, t \in (0,1),
\end{equation}
subject to the three-point BCs
\begin{equation}\label{eqE}
u'(0)=0,\; \beta u'(1) + u(\eta)=0,\; {\eta}\in [0,1],
\end{equation}
where $\beta >0$. The results of this Section are based on the manuscript~\cite{gijwnodea}.

One motivation for studying the BVP~\eqref{eqD}-\eqref{eqE} is that it occurs in some heat flow problems. This kind of problems were studied by Infante and Webb~\cite{gijwnodea}, who were motivated by earlier work by Guidotti and Merino~\cite{gm}. 

\begin{figure}[h]
\centering
\includegraphics[height=3cm,angle=0]{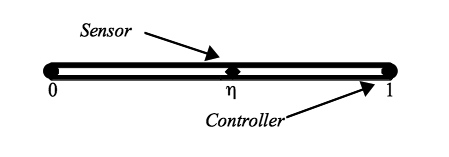}
\caption{A thermostat.}
\label{thermostat}
\end{figure}

In order to illustrate the physical interpretation of the BVP~\eqref{eqD}-\eqref{eqE}, suppose we have a heated bar of length 1. 
Then the temperature at a
point $x$ along the bar satisfies the one-dimensional heat equation
\begin{equation*}
u_t-u_{xx}=\hat{f}(t,x,u).
\end{equation*}
In the steady state, the equation becomes
\begin{equation*}
-u_{xx}=\hat{f}(x,u).
\end{equation*}
The use of the variable $t$ in lieu of the space variable $x$, gives
\begin{equation*}%\label{eqHF}
u''(t)+\hat{f}(t,u(t))=0,\ t \in (0,1).
\end{equation*}
The boundary conditions~\eqref{eqE}
can be interpreted as a model for a thermostat where in $t=0$ the bar is insulated and a controller at $t=1$ adds or removes heat
according to the temperatures detected by a sensor in
$t=\eta$ (see Figure \ref{thermostat}).
{}

In this simple model we have inserted along the bar only one sensor, but more complex models, with more controllers and sensors, may be studied. These kind of BCs are called~\emph{nonlocal} BCs and have received increasing attention in the last 20 years. As far as we know the
study of nonlocal BCs, in the context of ODEs, can be traced back
to Picone~\cite{Picone} in 1908, who considered multi-point BCs. For an
introduction to nonlocal problems we refer the reader to the
reviews~\cite{Conti, rma, sotiris, Stik, Whyburn}, the
papers~\cite{kttmna, ktejde, jw-gi-jlms} and the very well written notes~\cite{jw-seville}.

For further reading on thermostats problems with linear and nonlinear controllers, we refer the reader to 
 \cite{fama, df-gi-jp-prse, gijwems, gi-caa, kamkee, kamkee2, kapala, palamides, jwpomona, jwwcna04, webb-therm} and references therein.

In order to utilize the previous machinery, we construct the Green's function associated to the BVP~\eqref{eqD}-\eqref{eqE}, taking into account the presence of the nonlocal condition. 
Thus we consider the linear problem
\begin{align*}%\label{linBVP3}
u''(t)+y(t)=0,\ u'(0)=0,\; \beta u'(1) + u(\eta)=0.
\end{align*}
By integration, we obtain
$$u'(t)=u'(0)-\int_{0}^{t}y(s)ds$$ 
and, using the BC $u'(0)=0$  we get
$$u'(t)=-\int_{0}^{t}y(s)ds$$
and, by means of the Cauchy formula for iterated kernels, we obtain
$$u(t)=u(0)-\int_{0}^{t}\int_{0}^{w}y(s)dsdw=u(0)-\int_{0}^{t}(t-s)y(s)ds.$$
Therefore we have
$$u'(1)=-\int_{0}^{1}y(s)ds,$$
and
$$u(\eta)=u(0)-\int_{0}^{\eta}(\eta-s)y(s)ds.$$
Using the BCs, we have
$$0=\beta u'(1)+u(\eta)=-\beta \int_{0}^{1}y(s)ds+u(0)-\int_{0}^{\eta}(\eta-s)y(s)ds,$$
which, in turn, gives
$$
u(0)=\beta \int_{0}^{1}y(s)ds+\int_{0}^{\eta}(\eta-s)y(s)ds
$$
and
$$u(t)=\beta \int_{0}^{1}y(s)ds+\int_{0}^{\eta}(\eta-s)y(s)ds-\int_{0}^{t}(t-s)y(s)ds.$$
Thus we rewrite the BVP~\eqref{eqD}-\eqref{eqE} in the form~\eqref{eqhamm}, that is
$$
u(t)=\int_{0}^{1}k(t,s)f(u(s))\,ds :=Tu(t),
$$
where
\begin{equation}\label{ktherm}
k(t,s)=\beta +\begin{cases} \eta-s,\ &s\leq \eta,\\0,\ &s>\eta,
\end{cases}
-\begin{cases} t-s,\ &s\leq t,\\ 0,\ &s>t.
\end{cases}
\end{equation}

Here we discuss the case of $\beta+\eta\geq 1$ that leads to \emph{positive} solutions. We stress that a similar approach can be used to discuss the existence of solutions that \emph{change sign} (see for example \cite{gijwjiea} and \cite{gijwnodea}).

If $\beta +\eta \geq 1$ then $k(t,s)\geq 0$ for every $t,s \in
[0,1]$, and, since $k(t,s)$ is a decreasing function of $t$, we
have that the maximum of $k$ with respect to the variable $t$ is given by $k(0,s)$. Also the  the minimum of $k$ with respect to the variable $t$ is  $k(b,s)$ for $0\leq t \leq b$.
Thus we can take
\begin{equation*}
\Phi(s)=k(0,s)=\begin{cases}\beta,\ &s>\eta,\\ \beta+\eta - s,\
&s\leq \eta.
\end{cases}
\end{equation*}
If $\beta +\eta > 1$ we choose $[a,b]=[0,1]$, if $\beta +\eta = 1$
we choose $[a,b]=[0,b]$, with $b<1$. 

We have
\begin{equation*}
k(1,s)=\begin{cases}\beta-(1-s),\ &s>\eta,\\ \beta+\eta - 1,\
&s\leq \eta.
\end{cases}
\end{equation*}
For $\beta +\eta > 1$, we need  to choose $c$ so that $$ \beta + s
-1\geq c \beta, \text{ for }\eta<s\leq 1 $$ and $$ \beta + \eta
-1\geq c (\beta+\eta-s), \text{ for }0\leq s\leq \eta. $$

Hence it is sufficient to have $$c\leq 1-\frac{1}{\beta+\eta}. $$ For
$\beta +\eta = 1$, we have
\begin{equation*}
k(b,s)=\begin{cases}\beta +s -b,\ &s>\eta,\\ 1 - b,\ &s\leq \eta.
\end{cases}
\end{equation*}
Reasoning as in the previous case we see that it is enough to have
$$ c\leq 1-b. $$

The above calculations, in view of Theorem~\ref{eqhamm}, lead to an existence result for one or for
multiple solutions that are strictly positive on $[0,1)$. 

\begin{example}\label{extherm}
In the case $\beta +\eta > 1$, we can use  the cone
$$
K=\{u \in C[0,1], \min_{t \in [0,1]}u(t)\geq c \|u\|\},
$$
where $c= 1-\dfrac{1}{\beta+\eta}.$
A direct calculation gives
\begin{equation*}\label{38}
\frac{1}{m}=\beta+\eta^{2}/2\quad \text{and}\quad \frac{1}{M}
=(2\beta-1+\eta^{2})/2.
\end{equation*}
Take $\beta=1/2$ and $\eta=3/4$. This leads to $c=1/5$.

Then the $i_{K}(T,V_{\rho})=0$ condition needs
$$
f_{{\rho},{\rho / c}}\geq 32/9
$$
and
$i_{K}(T,K_{\rho})=1$
requires $$f^{0,\rho} \leq 32/25.$$
Therefore, provided that $f$ has a suitable growth, Theorem~\ref{eqhamm} can be applied.
\end{example}
\section{Nonlinear BCs}\label{subNBCs}
We now move to the case of nonlinear BCs and consider, as an illustrative example, a model of a chemical reactor. The results of this Section are based on the manuscript~\cite{genupa}.

The differential equation
\begin{equation}\label{eq1}
u''(t)-\lambda u'(t) + \lambda \mu (\beta - u(t)) e^{u(t)}=0,\ t \in (0,1),
\end{equation}
with the BCs
\begin{equation}\label{bcAhom}
u'(0)=\lambda u(0),\ u'(1)=0,
\end{equation} 
can be used as a model for the steady states of an adiabatic
chemical reactor of length 1. Here $\lambda$ is the Peclet
number, $\mu$ is the Damkohler number, $\beta$ is the
dimensionless adiabatic temperature rise and $u(t)$ is the local
temperature at a point $t$ of the tube, we refer the reader to
\cite{feng, gippreactor, Maam} and references therein.

Here we consider the more general BCs
\begin{equation}\label{bcA}
u'(0)=\lambda u(0),\ u'(1)=H[u],
\end{equation}
where $\beta, \lambda, \mu>0$ and $H$ is a suitable functional, not necessarily linear.

The nonlinear condition in~\eqref{bcA} can describe, for example,
a feedback control system on the reactor that adds or removes heat
according to the temperatures detected by some sensors located along the tube.

Due to the presence of the nonlinearity $H$, we seek solutions of the BVP~\eqref{eq1}-\eqref{bcA} by means of a \emph{perturbed} Hammerstein integral equation. This is quite a powerful trick that can be used in many situations, also when the BCs involve a \emph{linear} functionals, see for example~\cite{gi-caa, gijwems, jw-gi-jlms}. 

In our particular case, it is known that the solution of the BVP~\eqref{eq1}-\eqref{bcAhom} is given by 
$$
u(t)=\int_0^1 k(t,s) \mu (\beta - u(s)) e^{u(s)}\,ds, 
$$
where 
\begin{equation}\label{kerreact}
k(t,s)=\begin{cases}
e^{\lambda (t-s)},\ &s>t,\\
1,\ &s\leq t.
\end{cases}
\end{equation}
The Green's function~\eqref{kerreact} can be obtained by direct calculations and has been used in~\cite{feng, madbouly}. 

We seek the unique solution $\gamma$ of the linear BVP
\begin{equation*}%\label{lbvpaux}
\gamma''(t)-\lambda \gamma'(t)=0,\ \gamma'(0)=\lambda \gamma(0),\ \gamma'(1)=1,
\end{equation*}
which is given by 
$$\gamma(t)=\frac{1}{\lambda}e^{\lambda(t-1)}.$$
Therefore the solution of the BVP~\eqref{eq1}-\eqref{bcA} is given by the perturbed Hammerstein integral equation
\begin{equation}\label{phier}
u(t)={\gamma}(t)H[u] +\int_0^1 k(t,s) \mu (\beta - u(s)) e^{u(s)}\,ds.
\end{equation}

We prove the existence of strictly positive solutions (of norm less than $\beta$) of the integral equation~\eqref{phier} by solving, as we did in Section~\ref{fpin}, a slightly more general problem. In fact, we study equations of the form
\begin{equation}\label{eqpthamm}
u(t)={\gamma}(t)H[u] +\int_0^1 k(t,s)g(s)f(u(s))\,ds:=Tu(t).
\end{equation}
We make the following assumptions on the terms that occur in~\eqref{eqpthamm}.
\begin{itemize}
\item $f:[0, +\infty) \to [0, +\infty)$ is continuous.{}
\item  $k:[0,1] \times[0,1]\rightarrow (0,+\infty)$ is continuous.{}
\item There exist a continuous function $\Phi$ and a constant $c_{k} \in (0,1]$ such that
$$
c_{k}\Phi(s) \leq  k(t,s)\leq \Phi(s) \text{ for } t \in [0,1] \text{ and a.e. } \, s\in [0,1].
$${}
\item  $g \in L^1[0,1] $, $g(s) \geq 0$ for a.e. $s\in [0,1]$ and $\int_0^1 \Phi(s)g(s)\,ds >0$.{}
\item $\gamma \in C [0,1] $ and there exists $c_{\gamma} \in(0,1] \;\text{such that}\; \gamma(t) \geq c_{\gamma}\|\gamma\| \;\text{for}\; t \in [0,1]$.
\end{itemize}
Due to the hypotheses above, we are able to work in the cone
$$
K:=\{u\in C[0,1]:\ \min_{t\in [0,1]}u(t)\ge c\|u\|\},
$$
with $c=\min\{c_k,c_{\gamma}\}$ and we assume
\begin{itemize}
\item$H: K\to [0,+\infty)$ is compact.
\end{itemize}
Note that, since the range of $H$ is in $\mathbb{R}$, compact is the same as maps bounded sets to bounded sets (and continuous).

It is possible to show that the operator $T$ defined by~\eqref{eqpthamm} maps $K$ into $K$ and is compact.

In the following two Lemmas, rather than seeking \emph{global} linear bounds for the nonlinear functional $H$ we seek suitable \emph{local} linear bounds.

We begin with 
a condition which implies that the index is $1$.
\begin{lemma} \label{ind11}
Assume that
\begin{enumerate}
\item[$(\mathrm{I}^{1}_{\rho})$] there exist
 $\rho>0$, a linear functional $\alpha^{\rho}[\cdot]:K\rightarrow [0,+\infty)$ given by
$$
\alpha^{\rho}[u]=\int_0^1 u(t)\,dA^{\rho}(t)
$$
such that
\begin{itemize}
\item $dA^{\rho}$ is a positive Stieltjes measure with $A^{\rho}$ of bounded variation,
\item $\alpha^{\rho}[\gamma]<1$,
\item $H[u]\leq \alpha^{\rho}[u]$ for every $u\in \partial K_{
\rho}$,
\item the following inequality holds:
\begin{equation}\label{indice12}
  f^{c\rho,\rho}\Bigl(\sup_{t \in [0,1] }\Bigl\{\frac{\gamma(t)}{1-\alpha^{\rho}[\gamma]}\int_0^1\mathcal{K}^{\rho}(s)g(s)\,ds+\int_0^1k(t,s)g(s)\,ds\Bigr\} \Bigr)< 1,
\end{equation}
where
\begin{equation*}
 f^{c\rho,\rho}:=\sup\Bigl\{\frac{f(u)}{\rho},\,\, c\rho\le
u \le \rho\Bigr\}
\,\text{  and    } \,\mathcal{K}^{\rho}(s):=\int_0^1 k(t,s)\,dA^{\rho}(t).
\end{equation*}
\end{itemize}
\end{enumerate}
Then $i_{K}(T,K_{\rho})$ is $1$.
\end{lemma}

\begin{proof}
Note that if
$u \in \partial K_{\rho}$  then we have $c\rho\le u(t) \le \rho$ for every $t \in [0,1]$.

We show that $\mu\, u\neq Tu$ for every $u\in \partial K_{\rho}$ and for every $\mu \geq 1$; this ensures that the index is 1 on $K_{\rho}$. In fact, if this does
not happen, there exist $\mu \geq 1$ and $u\in \partial K_{\rho}$ such that, for every $t \in [0,1]$,
\begin{equation*}
\mu\, u(t)= Tu(t)=\gamma (t)H[u]+\int_{0}^{1}k(t,s)g(s)f(u(s))ds.
\end{equation*}
Then we have
\begin{equation}\label{rel}
\mu\, u(t)\le \gamma(t)\alpha^{\rho}[u]+\int_{0}^{1} k(t,s)g(s)f(u(s))ds.
 \end{equation}
Applying $\alpha^{\rho}$ to the both sides of \eqref{rel} gives
\begin{equation*}
\mu\, \alpha^{\rho}[u]\le \alpha^{\rho}[\gamma]\alpha^{\rho}[u]+\int_0
^1\mathcal{K}^{\rho}(s)g(s)f(u(s))ds.
 \end{equation*}
 Thus we have
\begin{equation} \label{alpha}
\alpha^{\rho}[u]\le
\frac{1}{\mu-\alpha^{\rho}[\gamma]}\int_0
^1\mathcal{K}^{\rho}(s)g(s)f(u(s))ds\le\frac{1}{1-\alpha^{\rho}[\gamma]}\int_0
^1\mathcal{K}^{\rho}(s)g(s)f(u(s))ds.
 \end{equation}
 Using \eqref{alpha} in \eqref{rel} we obtain
\begin{align*}
\mu\, u(t)&
\le\frac{\gamma(t)}{1-\alpha^{\rho}[\gamma]}\int_0
^1\mathcal{K}^{\rho}(s)g(s)f(u(s))ds+\int_{0}^{1}
k(t,s)g(s)f(u(s))ds\\
&\le\rho f^{c\rho,\rho}\Bigl(\frac{\gamma(t)}{1-\alpha^{\rho}[\gamma]}\int_0^1\mathcal{K}^{\rho}(s)g(s)ds+\int_{0}^{1} k(t,s)g(s)ds\Bigr).
 \end{align*}
 Taking the supremum in $[0,1]$ gives
\begin{align*}
\mu\rho&\le
\rho f^{c\rho,\rho}\Bigl(\sup_{t \in[0,1]}\Bigl\{\frac{\gamma(t)}{1-\alpha^{\rho}[\gamma]}\int_0
^1\mathcal{K}^{\rho}(s)g(s)ds+\int_{0}^{1}
k(t,s)g(s)ds\Bigr\}\Bigr)
\end{align*}
and using the hypothesis \eqref{indice12} we can conclude that
$\mu \rho <\rho$. This contradicts the fact that $\mu \geq 1$ and proves the result.
\end{proof}
Now we give a  condition which implies that the index is $0$ on the set $V_{\rho}$.
\begin{lemma} \label{indice0*}
Assume that
\begin{enumerate}
\item[$(\mathrm{I}^{0}_{\rho})$] there exist
 $\rho>0$, a linear functional $\alpha^{\rho}[\cdot]:K\rightarrow [0,+\infty)$ given by
 $$\alpha^{\rho}[u]=\int_0^1 u(t)\,dA^{\rho}(t)$$
such that
\begin{itemize}
\item $dA^{\rho}$ is a positive Stieltjes measure with $A^{\rho}$ of bounded variation,
\item $\alpha^{\rho}[\gamma]<1$,
\item $H[u]\geq \alpha^{\rho}[u]$ for every $u\in \partial V_{
\rho}$,
\item the following inequality holds:
\begin{equation}\label{ind 0 2 +}
  f_{\rho,\rho/c}\Bigl(\inf_{t \in [0,1] }\Bigl\{\frac{\gamma(t)}{1-\alpha^{\rho}[\gamma]}\int_0^1\mathcal{K}^{\rho}(s)g(s)\,ds+\int_0^1 k(t,s) g(s)\,ds\Bigr\} \Bigr)> 1.
\end{equation}
\end{itemize}
\end{enumerate}
Then $i_{K}(T,V_{\rho})$ is $0$.
\end{lemma}

\begin{proof}
Note that the constant function $e(t)\equiv
1$  for $t\in[0,1]$ belongs to $K$. Furthermore observe that if $u \in \partial V_{\rho}$ then we have $\rho\le u(t) \le \rho/c$ for
every $t \in [0,1]$.

We prove that $u\not=Tu+\lambda e$ for every $u\in \partial V_{\rho}$ and for every $\lambda \geq 0$; this ensures that the index is $0$ on $V_{\rho}$.

Let $u\in \partial V_{\rho}$ and $\lambda \geq 0$ such that
$u=Tu+\lambda\,e$. Then we have, for $t\in[0,1]$,
\begin{align}\label{Eq0}
   u(t)=&\gamma(t)H[u]+\int_{0}^{1} k(t,s)g(s)f(u(s))ds+ \lambda e(t)\\
 \nonumber \geq & \gamma(t)\alpha^{\rho}[u]+\int_{0}^{1}k(t,s)g(s)f(u(s))ds.
  \end{align}
Thus we have
\begin{equation*}
\alpha^{\rho}[u]\geq \alpha^{\rho}[\gamma]\alpha^{\rho}[u]+\int_{0}^{1}
\mathcal{K}^{\rho}(s)g(s)f(u(s))\,ds.
\end{equation*}
This implies
\begin{equation}\label{Eq1}
  \alpha^{\rho} [u]\geq  \frac{1}{1-\alpha^{\rho}[\gamma] }\int_{0}^{1} \mathcal{K}^{\rho}(s)g(s)f(u(s))\,ds.
\end{equation}
Using \eqref{Eq1} in \eqref{Eq0}  we obtain
\begin{align*}
 u(t)&\geq
 \gamma(t)\frac{1}{1-\alpha^{\rho}[\gamma] }\int_{0}^{1} \mathcal{K}^{\rho}(s)g(s)f(u(s))\,ds+\int_{0}^{1} k(t,s) g(s)f(u(s))\,ds\\
& \geq \rho f_{\rho,\rho/c}\Bigl(\frac{\gamma(t)}{1-\alpha^{\rho}[\gamma]}\int_0^1\mathcal{K}^{\rho}(s)g(s)\,ds+\int_0^1
k(t,s)g(s)\,ds\Bigr).
 \end{align*}
Taking the infimum for $t\in [0,1]$ gives
\begin{equation*}
\rho \geq \rho f_{\rho,\rho/c} \Bigl(\inf_{t \in [0,1]
}\Bigl\{\frac{\gamma(t)}{1-\alpha^{\rho}[\gamma]}\int_0^1\mathcal{K}^{\rho}(s)g(s)\,ds+\int_0^1
k(t,s) g(s)\,ds\Bigr\} \Bigr).
 \end{equation*}
  Thus from \eqref{ind 0 2 +}  we have $\rho> \rho$.
This is a contradiction that proves the result. 
\end{proof}
\begin{remark}
When $H[u]\equiv 0$, the growth condition~\eqref{indice12} reads more simply 
\begin{equation*}
  f^{c\rho,\rho}< m,
\end{equation*}
where
$$
\frac{1}{m}=\sup_{t\in [0,1]} \int_{0}^{1} k(t,s)g(s) \,ds.
$$
while the growth condition~\eqref{ind 0 2 +} reads 
\begin{equation*}
  f_{\rho,\rho/c}> M.
\end{equation*}
where
$$
\frac{1}{M} =\inf_{t\in [0,1]}\int_{0}^{1}k(t,s)g(s)\,ds.
$$
Note that this setting also improves the result in Example~\ref{extherm}, allowing more freedom with the choice of the nonlinearity.
\end{remark}
A Theorem similar to Theorem~\ref{thmmsol1} holds for the integral equation \eqref{eqpthamm}, yielding existence of \emph{strictly positive} solutions, we omit the statement of this result.

We turn our attention back to the BVP  \eqref{eq1}-\eqref{bcA} and we seek solutions of norm less than~$\beta$, by studying the integral equation
\begin{equation*}%\label{loc-ham}
u(t)=\gamma(t)H[u]+\int_{0}^{1}k(t,s)f(u(s))\,ds.
\end{equation*}
where
\begin{equation*}
\gamma(t)=\frac{1}{\lambda}e^{\lambda(t-1)},\quad
k(t,s)=\begin{cases}
e^{\lambda (t-s)},\ &s>t,\\
1,\ &s\leq t,
\end{cases}
\end{equation*}
and
\begin{equation*}
f(u)=\begin{cases}
\mu (\beta - u) e^{u},\ &u\leq \beta, \\
0,\ & u>\beta.
\end{cases}
\end{equation*}
We work in the cone
\begin{equation*}%\label{eqcone}
K = \{u\in C[0,1], \; \min_{t\in [0,1]}u(t) \geq c \|u\|\},
\end{equation*}
 where the constant $c= e^{-\lambda}$, since
$$
e^{-\lambda} \leq k(t,s)\leq 1\text{ for } t \in [0,1]\times [0,1],
$$
and the conditions on $k$ and $\gamma$ are satisfied with  $\Phi(s)=1$ and $c_k=c_{\gamma}=e^{-\lambda}$.
\begin{example}
In order to illustrate the growth conditions, we consider the BVP
\begin{equation}\label{ex-chem-eq}
u''(t)-\frac{1}{3} u'(t) + \frac{3}{10} \Bigl(\frac{11}{5} - u(t)\Bigr) e^{u(t)}=0,\ t \in (0,1),
\end{equation}
\begin{equation}\label{ex-chem-bc}
u'(0)=\frac{1}{3} u(0),\, u'(1)=10^{-\frac{3}{2}}\sqrt{u(1/2)}.
\end{equation}

The choice $$\rho_1=\frac{71}{1000},
\,\,\rho_2=\frac{53}{25},
\,\,\alpha^{\rho_1}[u]=\frac{1}{10}u(1/2)
,\,\,\alpha^{\rho_2}[u]=10^{-\frac{5}{4}}u(1/2),$$
yields (in what follows the numbers are rounded to the third decimal place unless exact)
\begin{itemize}
\item[] $\alpha^{\rho_1}[\gamma]=0.254<1$ and $\alpha^{\rho_2}[\gamma]=0.143<1$,
\item[]  $H[u]=10^{-\frac{3}{2}}\sqrt{u(1/2)}
\geq \frac{1}{10}u(1/2)=\alpha^{\rho_1}[u]$ for $\rho_1\leq
u\leq\rho_1/c$,
\item[]  $H[u]=10^{-\frac{3}{2}}\sqrt{u(1/2)}\leq 10^{-\frac{5}{4}}u(1/2)=\alpha^{\rho_2}[u]$ for $c\rho_2\leq u\leq\rho_2$,
\item[]  $\inf\left\{f(u):\,\,u\in [\rho_1,
\rho_1/c]\right\}=2.057>\frac{71}{1000}\cdot 1.917$,
\item[]  $\sup\left\{f(u):\,\,u\in \left[c\rho_2,
\rho_2\right]\right\}=2.811<\frac{53}{25}\cdot 2.551$.
\end{itemize}
Thus the conditions $(\mathrm{I}^{0}_{\rho_1})$ of Lemma
\ref{indice0*} and $(\mathrm{I}^{1}_{\rho_2})$ of Lemma
\ref{ind11} are satisfied. Then it follows that
the BVP~\eqref{ex-chem-eq}-\eqref{ex-chem-bc} has a strictly positive solution $u \in {K}_{\rho_2}\setminus \overline{V}_{\rho_1}$ with the following localization property:
$$
\rho_1=71/1000 \leq u(t)\leq 53/25= \rho_2,\ \text{for every}\
t\in [0,1].
$$
\end{example}
\chapter{Radial solutions of PDEs}\label{radsolpdes}
We now briefly illustrate how to apply the previously theory in order to deal with the existence of radial solutions of systems of elliptic PDEs. In particular we study the case of annular and exterior domains; a reader interested in this topic might find interesting the review~\cite{jj-ks} and the papers~\cite{genupa2, dolo, do6, gipp-nodea, kqljwjde, webb}.

The methodology here is to associate to the elliptic system a system of Hammerstein integral equations of the type
\begin{gather}
\begin{aligned}\label{systham}
u(t)=&\int_{0}^{1}k_1(t,s)g_1(s)f_1(u(s),v(s))\,ds, \\
v(t)=&\int_{0}^{1}k_2(t,s)g_2(s)f_2(u(s),v(s))\,ds,%
\end{aligned}
\end{gather}
a form a little more general than~\eqref{ham-sys}. 

We make the following assumptions on the terms that occur in~\eqref{systham}, for $i=1,2$.
\begin{itemize}
\item $f_i:  [0, +\infty)\times [0, +\infty) \to [0, +\infty)$ is continuous.
\item  $k_i:[0,1] \times[0,1]\rightarrow [0,+\infty)$ is continuous.{}
\item There exist a continuous function  $\Phi_i:[0,1]\to
[0, +\infty)$, an interval $[a_i,b_i]\subset [0,1]$ and a constant $c_i \in (0,1]$ such that{}
\begin{align*}
  k_i(t,s)\leq \Phi_i(s) \text{ for }& t \in [0,1] \text{ and a.e. }s \in
  [0,1],\\
  c_i \Phi_i(s)\leq k_i(t,s) \text{ for }& t \in [a,b] \text{ and  a.e. }s \in
  [0,1].
\end{align*}{}
\item  $g_i \in L^1[0,1] $, $g_i(s) \geq 0$ for a.e. $s\in [0,1]$ and $\int_{a_{i}}^{b_{i}} \Phi_i(s)g_i(s)\,ds >0$.{}
\end{itemize}

Under the assumptions above, we may proceed in a similar way as in Section~\ref{secsyst} and look for solutions of the system~\eqref{systham} in the cone
\begin{equation*}
 K:=  \{ (u,v) \in \tilde{K_1} \times \tilde{K_2}  \},
\end{equation*}
where
$$
\tilde{K_i}:= \{ w \in C[0,1]: w(t)\geq 0\ \text{and}\; \min_{t\in [a_i,b_i]} w(t) \geq c_i
\|w\| \}.
$$
Results similar to Lemmas \ref{ind1b-sys}, \ref{idx0b1-sys}, \ref{idx0b3-sys} and Theorem~\ref{mult-sys} hold in this context. For brevity, we do not state these results and refer to~\cite{gipp-nodea}, but, nevertheless, we point out that the main difference lies within the constants involved, that take into account (in a similar way as in Section~\ref{subNBCs}) the term $g_i$, namely
$$
\frac{1}{m_i}=\sup_{t\in [0,1]} \int_{0}^{1} k_i(t,s)g_i(s) \,ds,\quad   \frac{1}{M_i}=\inf_{t\in
[a_i,b_i]}\int_{a_i}^{b_i} k_i(t,s) g_i(s) \,ds.
$$
\section{Radial solutions of systems in annular domains}\label{radann}
Consider the systems of BVPs
\begin{gather}
\begin{aligned}\label{ellbvp-secapp}
\Delta u + h_1(|x|) f_1(u,v)=0,\ |&x|\in [R_1,R_0], \\
\Delta v + h_2(|x|) f_2(u,v)=0,\ |&x|\in [R_1,R_0],\\
\frac{\partial u}{\partial r}\Bigr\rvert_{\partial
B_{R_0}}=0\ \text{and}\
(u(R_1 x)-\beta  &u(R_\eta x))\Big|_{x\in\partial B_{1}}=0,\\
v\Bigr\rvert_{\partial
B_{R_0}}=0 \ \text{and}\
 & \frac{\partial v}{\partial r}\Big|_{\partial B_{R_1}}=0,%
\end{aligned}
\end{gather}
where $x\in \mathbb{R}^n $, $n\geq 2$,
$0<R_1<R_0<+\infty$, $R_\eta \in (R_1,R_0)$ and $0\leq \beta <1$ and  $\dfrac{\partial}{\partial r}$ denotes (as in~\cite{nirenberg}) differentiation in the radial direction $r=|x|$.

We assume that for $i=1,2$,
\begin{itemize}
\item $f_i:  [0, +\infty)\times [0, +\infty) \to [0, +\infty)$ is continuous.
{}
\item $h_i:  [R_1,R_0] \to [0, +\infty)$ is continuous.
\end{itemize}

In order to deal with the system~\eqref{ellbvp-secapp}, consider  in $\mathbb{R}^n$, $n\ge 2$, the equation
\begin{equation}\label{eqell}
\triangle w+h(|x|)f(w) = 0, \  \text{for a.e.}\  |x|\in
[R_{1},R_{0}].
\end{equation}

To establish the existence of radial solutions $w=w(r)$,
$r=|x|$, we proceed as in~\cite{kqljlms, lan-lin-na, kqljwjde} and
rewrite \eqref{eqell} in the form
\begin{equation}\label{eqinterm}
w''(r) + \dfrac{n-1}{r}w'(r) + h(r)f(w(r))= 0
\quad\text{a.e. on } [R_{1}, R_{0}].
\end{equation}
Set
$w(t)=w(r(t))$, where, for $t\in[0,1]$,
\begin{equation*}
r(t):=\begin{cases}
R_0^{1-t}R_1^{t},\,\,\,\,\,\,\,\,\,\,\,\,\,\,\,\,\,\, n=2,\\
({R_{0}^{-(n-2)}}+({R_{1}^{-(n-2)}}-{R_{0}^{-(n-2)}})t)^{-1/(n-2)},\ &n\geq 3.
\end{cases}
\end{equation*}
Take, for $t\in[0,1]$,
\begin{equation*}
\phi(t):=\begin{cases}
r^2(t) \log^2(R_0/R_1),\,\,\,\,\,\,\,\,\,\,\,\,\,\,\,\,\,\, n=2,\\
\Bigl(\frac{R_{1}^{-(n-2)}-R_{0}^{-(n-2)}}{n-2}\Bigr)^2 \Bigl(R_{0}^{-(n-2)}+(R_{1}^{-(n-2)}-R_{0}^{-(n-2)})t \Bigr)^{\frac{-2(n-1)}{n-2}},\ n\geq 3,
\end{cases}
\end{equation*}
then~\eqref{eqinterm}
becomes
\begin{equation*}%\label{eqdef}
w''(t) + {\phi}(t) h(r(t)) f(w(t)) = 0, \  \text{a.e. on}\
[0,1].
\end{equation*}
 Set
$u(t)=u(r(t))$ and
$v(t)=v(r(t))$. Thus, to the system \eqref{ellbvp-secapp} we associate the system of ODEs
\begin{gather}
\begin{aligned}\label{1syst}
u''(t) + g_1(t) f_1(u(t),v(t)) = 0, \quad \text{a.e. on } [0,1], \\
v''(t) + g_2(t) f_2(u(t),v(t)) = 0, \quad \text{a.e. on } [0,1],%
\end{aligned}
\end{gather}
subject to the BCs
\begin{gather}
\begin{aligned}\label{1BC}
u'(0)=0,\ u(1)&={\beta}u({\eta}), \\
v(0)=0,\ v'(1)&=0,
\end{aligned}
\end{gather}
where
\begin{equation*}
g_i(t):={\phi}(t) h_i(r(t)),
\end{equation*}
and $\eta \in (0,1)$ is such that $r(\eta)=R_{\eta}$.

Therefore we can study the existence of radial solutions of the
system~\eqref{1syst}-\eqref{1BC} by means of the system~\eqref{systham}
where $k_1$ is given by
\begin{equation*} %\label{ker1}
k_1(t,s)=\dfrac{1}{1-\beta}(1-s)-\begin{cases}
\dfrac{\beta}{1-\beta}(\eta -s), &  s \le \eta,\\ \quad 0,&
s>\eta,
\end{cases}
 - \begin{cases} t-s, &s\le t, \\ \quad 0,&s>t,
\end{cases}
\end{equation*}
and $k_2$ is given by~\eqref{kRob}. 

Note that the kernel $k_1$ is non-negative when $0\leq \beta <1$.
Upper and lower bounds for $k_1$ were carefully studied in~\cite{jw-wilm}, where it was shown that 
one may use as $[a_1,b_1]=[0,\eta],$
$$ 
\Phi_1(s)=\begin{cases}\quad \quad \dfrac{1-s}{1-\beta},
& \text{ if } \eta <s \leq 1,\\
\dfrac{1-s-\beta(\eta-s)}{1-\beta},
& \text { if } 0 \leq s \leq \eta.
\end{cases}
$$
and 
$$c_1=\frac{(1-\eta)}{(1-\beta\eta)}.$$
\section{Radial solutions in exterior domains}
We now consider the systems of BVPs
\begin{gather}
\begin{aligned}\label{ellbvpsec}
\Delta u + h_1(|x|) f_1(u,v)=0,&\ |x|\in [R_1,+\infty), \\
\Delta v + h_2(|x|) f_2(u,v)=0,&\ |x|\in [R_1,+\infty),\\
u(R_1x)=\beta u(R_\eta x)\,\,\text{for  }x\in \partial &B_{1},\  \displaystyle\lim_{|x|\rightarrow +\infty} u(|x|)=0,\\
v\Bigr\rvert_{\partial
B_{R_1}}=0,&\ \displaystyle\lim_{|x|\to+\infty}v(|x|)=0,%
\end{aligned}
\end{gather}
where $x\in \mathbb{R}^n$, $n\geq 3$, $\beta\geq 0$, $R_1>0$, $R_\eta, R_\xi \in (R_1,+\infty)$.

We assume that the following holds, for $i=1,2$.
\begin{itemize}
\item $f_i:  [0, +\infty)\times [0, +\infty) \to [0, +\infty)$ is continuous.
\item  $h_i:[R_1,+\infty) \to [0,+\infty)$ is continuous  and $h_i(|x|)\leq \frac{1}{|x|^{n+\mu_i}}$ for $|x|\to +\infty$ and for some $\mu_i>0$.
{}
\end{itemize}

In a similar way as in Section~\ref{radann} we consider  in $\mathbb{R}^n$, $n\ge 3$, the equation
\begin{equation}\label{eqell2}
\triangle w+ h(|x|)f(w) = 0, \  \text{for a.e.}\  |x|\in
[R_{1},+\infty).
\end{equation}

In order to establish the existence of radial solutions $w=w(r)$, $r=|x|$, we proceed as in~\cite{but} and we  rewrite \eqref{eqell2} in the form
\begin{equation}\label{eqinterm2}
w''(r) + \dfrac{n-1}{r}w'(r) + h(r)f(w(r))= 0,  \text{  for }\  r\in [R_{1}, +\infty).
\end{equation}
Set $w(t)=w(r(t))$, where, for $t\in[0,1]$,
\begin{equation*}
r(t):=R_1\,t^{\frac{1}{2-n}}.
\end{equation*}
Take, for $t\in[0,1]$,
\begin{equation*}
\phi(t):=r(t)\,\frac{R_1}{(n-2)^2}\,t^{\frac{2n-3}{2-n}},
\end{equation*}
then the equation~\eqref{eqinterm} becomes
\begin{equation*}
w''(t) + {\phi}(t) h(r(t)) f(w(t)) = 0, \  \text{  on}\ [0,1].
\end{equation*}
Set $u(t)=u(r(t))$ and $v(t)=v(r(t))$. 

Thus, to the system \eqref{ellbvpsec} we associate the system of ODEs
\begin{gather}
\begin{aligned}\label{5syst}
u''(t) + g_1(t) f_1(u(t),v(t)) =& 0,\  t\in (0,1), \\
v''(t) + g_2(t) f_2(u(t),v(t)) =& 0, \ t\in  (0,1),\\
\end{aligned}
\end{gather}
with BCs
\begin{gather}
\begin{aligned}\label{5systBC}
u(0)=0,\ u(1)={\beta}u({\eta})&, \; 0 < {\eta}<1, \\
v(0)=v(1)=&0,
\end{aligned}
\end{gather}
where
$$
g_i(t):={\phi}(t) h_i(r(t)),
$$
and $\eta \in (0,1)$ is such that $r(\eta)=R_{\eta}$.

We study the existence of solutions of the system~\eqref{5syst}-\eqref{5systBC} via the Hammerstein integral system~\eqref{systham},
where, this time, 
\begin{equation*}
k_1(t,s)=\dfrac{1}{1-\beta \eta} t(1-s)-
\begin{cases}
\dfrac{\beta t}{1-\beta \eta}(\eta -s), &  s \le \eta\\ 
\quad 0,& s>\eta
\end{cases}
 - \begin{cases}
  t-s, &s\le t, \\ 
  \quad 0,&s>t,
\end{cases}
\end{equation*}
and $k_2$ is given by~\eqref{kDir}.

Note that the kernel $k_1$ is non-negative for $0\leq \beta \eta <1$. A careful study of the upper and lower bounds for $k_1$ was done, once again, in~\cite{jw-wilm}. These results can be summarized as follows.

When $\beta\leq 1$, one may use
$$
\Phi_1(s)=
\begin{cases}
\dfrac{1}{1-\beta\eta} s(1-s), &\text{ if } s>\eta,\\
\dfrac{1}{1-\beta\eta} s\bigl(1-s-\beta(\eta-s)\bigr), &\text{
if } 0 \leq s \leq \eta.
\end{cases}
$$
and
$c_1=\min\Bigl\{\dfrac{1-\eta}{1-\beta\eta},\beta\eta\Bigr\}$.

When $1<\beta < 1/\eta$, one may take
$$
\Phi_1(s)=
\begin{cases}
\dfrac{1}{1-\beta\eta} s(1-s), &\text{ if } s>\beta\eta,\\
\dfrac{1}{1-\beta\eta} \beta\eta(1-s), &\text{ if } \eta < s <
\beta\eta,\\ \dfrac{1}{1-\beta\eta} \beta(1-\eta)s, &\text{ if
} 0 \leq  s \leq \eta,
\end{cases}
$$
and $c_1=\eta$.
\section*{Conclusions and further reading}\addcontentsline{toc}{section}{Conclusions and further reading}
We have briefly shown that, in some cases, the existence of radial, non-negative solutions of systems of elliptic PDEs subject to local and nonlocal BCs, can be studied via systems of Hammerstein integral equations. Therefore, provided that the nonlinearities involved have a suitable growth, existence, multiplicity and non-existence results can be obtained. 
Finally, we mention that it is possible to tailor this theory in order to deal, in the spirit of Section~\ref{subNBCs}, with elliptic systems with more general nonlinear BCs, we refer the reader to the papers~\cite{genupa, genupa2}.

\chapter*{Acknowledgments}\addcontentsline{toc}{chapter}{Acknowledgments}

G. Infante would like to thank the Departamento de An\'alise Matem\'atica of the Universidade de Santiago de Compostela, 
the Department of Mathematical Analysis of the University of Ruse and J. A. Cid, R. Figueroa and F. A. F. Tojo (Organizers of the Workshop ``Differential Equations and Applications'') for their warm hospitality, generous support and the opportunity to deliver these notes.
G. Infante would also like to thank F. A. F. Tojo, P. Pietramala and J. R. L. Webb, for carefully checking some drafts of these notes.
G. Infante was partially supported by G.N.A.M.P.A. - INdAM (Italy) and the Erasmus$+$ program.

\end{document}